\documentclass[12pt]{amsart}

\usepackage{amssymb,epic,eepic}

\usepackage[cmtip,matrix,arrow]{xy}

\newtheorem{theorem}{Theorem}[section]

\newtheorem{proposition}[theorem]{Proposition}

\newtheorem{lemma}[theorem]{Lemma}

\theoremstyle{remark}
\newtheorem{remark}[theorem]{Remark}
\newtheorem{definition}[theorem]{Definition}

\newtheorem{example}[theorem]{Example}

\newcommand\A{\mathcal{A}}
\newcommand\be{\begin{equation}}
\newcommand\ee{\end{equation}}
\newcommand\M{\mathcal{M}}
\renewcommand\L{\mathcal{L}}

\renewcommand{\O}{\mathcal{O}}

\newcommand{\JJ}{{J}}

\newcommand{\U}{\on{U}}

\newcommand{\R}{\mathbb{R}}
\newcommand{\C}{\mathbb{C}}

\newcommand{\Z}{\mathbb{Z}}
\newcommand{\Q}{\mathbb{Q}}

\newcommand{\nh}{{}}

\newcommand\lie[1]{\mathfrak{#1}}

\newcommand{\g}{\lie{g}}

\renewcommand{\t}{\lie{t}}
\newcommand{\Alc}{\lie{A}}

\newcommand{\on}{\operatorname}

\newcommand{\curv}{ \on{curv} } 
\newcommand{\Ad}{ \on{Ad} } \newcommand{\Hol}{ \on{Hol} }
\newcommand{\Eul}{ \on{Eul} }

\renewcommand{\ker}{ \on{ker}} 
 \newcommand{\Spin}{ \on{Spin}}
\newcommand{\SU}{ \on{SU}}

\newcommand{\Oma}{{\underline \Om}}
\renewcommand{\index}{ \on{index}}
\newcommand{\codim}{\on{codim}}

\newcommand{\D}{ \mathcal{D} }

\newcommand\dirac{/\kern-1.2ex\partial} 
\newcommand\qu{/\kern-.7ex/} 
\newcommand{\Waff}{W_{\on{aff}}} 

\newcommand{\lra}{\longrightarrow}
\newcommand{\hra}{\hookrightarrow}

\renewcommand{\d}{{\mbox{d}}}
\newcommand{\ol}{\overline}

\newcommand\Phinv{\Phi^{-1}}

\newcommand\lam{\lambda}
\newcommand\Lam{\Lambda}
\newcommand\Sig{\Sigma}
\newcommand\sig{\sigma}
\newcommand\eps{\epsilon}
\newcommand\Om{\Omega}
\newcommand\om{\omega}

\newcommand{\f}{\frac}

\renewcommand{\H}{\ca{H}}
\newcommand{\p}{\partial}
\renewcommand{\l}{\langle}
\renewcommand{\r}{\rangle}
\newcommand{\hh}{{\textstyle \f{1}{2}}}
\newcommand{\ti}{\tilde}
\newcommand{\olt}{\overline{\theta}}

\newcommand\Td{\on{Td}}
\newcommand\Ch{\on{Ch}}

\newcommand\beqn{\begin{equation}}      
\newcommand\eeqn{\end{equation}}      
\newcommand{\ca}{\mathcal}
\newcommand{\wh}{\widehat}
\newcommand{\wt}{\widetilde}
\newcommand{\mf}{\mathfrak}
\newcommand{\beq}{\begin{eqnarray*}}
\newcommand{\eeq}{\end{eqnarray*}}
\begin{document}

\title{A fixed point formula for loop group actions}

\date{\today}

\author{A. Alekseev}
\address{Institute for Theoretical Physics \\ Uppsala University \\
Box 803 \\ \mbox{S-75108} Uppsala \\ Sweden}
\email{alekseev@teorfys.uu.se}

\author{E. Meinrenken}
\address{University of Toronto, Department of Mathematics,
100 St George Street, Toronto, Ontario M5R3G3, Canada }
\email{mein@math.toronto.edu}

\author{C. Woodward}
\address{Mathematics-Hill Center, Rutgers University,
110 Frelinghuysen Road, Piscataway NJ 08854-8019, USA}
\email{ctw@math.rutgers.edu}

\begin{abstract}  We express the index of the Dirac operator on 
symplectic quotients of a Hamiltonian loop group manifold with proper
moment map in terms of fixed point data.
\end{abstract}

\maketitle


\section{Introduction}
\label{sec:intro}

The purpose of this paper is to generalize to loop group manifolds the
following result on symplectic quotients of Hamiltonian actions of
compact groups.  Consider a compact symplectic manifold $M$, with a
Hamiltonian action of a compact 
%
connected 
Lie group $G$ and equivariant
pre-quantum line bundle. Associated to these data is a (virtual)
character $\chi(M)$ of $G$, defined as the equivariant index of a
$\Spin_c$-Dirac operator. $\chi(M)$ can be computed from the  
Atiyah-Segal-Singer theorem in terms of fixed point data, or 
from the {\em quantization commutes with reduction} principle
in terms of indices of symplectic quotients.  A combination of these two
expressions leads to formulas for indices of symplectic quotients as a
sum of fixed point contributions.

This paper is concerned with a similar fixed point formula for
symplectic quotients of pre-quantized Hamiltonian loop group 
manifolds $\wh{M}$ with proper moment map. While $\wh{M}$ 
itself is infinite-dimensional, the properness assumption 
implies that its symplectic quotients are 
compact. We will not attempt to make sense of the equivariant index in 
infinite dimensions, or to define fixed point contributions on 
$\wh{M}$. Instead, we consider a finite dimensional
compact $G$-manifold $M$, obtained from $\wh{M}$ as a quotient by the
{\em based} loop group $\Om G\subset LG$. 
Our main result (Theorem \ref{onal}) is a formula for 
indices of symplectic quotients of $\wh{M}$ in terms of fixed 
point data on $M$. 
The fixed point contributions are reminiscent 
of the right hand side of the equivariant index theorem. 
However, $M$ does not carry a naturally induced 
$\Spin_c$-structure, since neither the symplectic 
form nor the line bundle descend to $M$. 
Heuristically, the formula follows by application of 
the equivariant index theorem to the loop group manifold, 
and subsequent ``renormalization'' of the
infinities on both sides. 
In a companion paper we use Theorem \ref{onal} to  
compute Verlinde numbers for moduli spaces of flat 
connections on surfaces.

A project with similar goals was undertaken earlier by S. Chang.  His
unfinished manuscript, dealing with the special case that the moment
map is transversal to the Cartan subalgebra, may be found at
\cite{ch:fi}. Chang's formula gives the indices of symplectic
quotients as a fixed point formula for a torus action on a compact
symplectic space, the ``imploded cross-section'' of $\wh{M}$.

The contents of the paper are as follows. In Section \ref{sec:spinc}
we discuss $\Spin_c$-quantization of finite dimensional Hamiltonian
manifolds, and describe the finite dimensional version of our fixed
point formula. Section \ref{sec:review} is dedicated to a review of
loop group actions and group-valued moment maps. In Section
\ref{sec:mainresult} we state the main theorem and in Section
\ref{sec:theproof} we describe the proof.
\vskip.6in
\begin{center}
{\bf Notation}
\end{center}
\vskip.1in

Throughout the paper $G$ will denote a compact, connected Lie group,
and $\g$ its Lie algebra. We denote by $R(G)$ the ring of characters 
of finite-dimensional representations.
We let $T$ be a maximal torus in $G$, and $\t$ its
Lie algebra. The integral lattice $\Lambda\subset \t$ is defined as
the kernel of the exponential map $\exp:\t\to T$, and the (real)
weight lattice $\Lambda^*\subset\t^*$ is its dual.  Embed
$\t^*\hra\g^*$ as the fixed point set for the coadjoint action of
$T$.  Every $\mu\in\Lambda^*$ defines a 1-dimensional
$T$-representation, denoted $\C_\mu$, where $t=\exp\xi$ acts by
$t^\mu:=e^{2\pi i \l\mu,\xi\r}$.  This representation extends uniquely
to the stabilizer group $G_\mu$ of $\mu$.  We let $W$ be the Weyl
group of $(G,T)$ and $\mf{R}\subset\Lambda^*$ the set of roots. We fix
a set of positive roots $\mf{R}_+\subset\mf{R}$ and let
$\t_+\subset\t$ and $\t^*_+\subset\t^*$ be the corresponding positive
Weyl chambers. For any dominant weight
$\mu\in\Lambda^*_+:=\Lambda^*\cap\t^*_+$ we denote by $V_\mu$ the
irreducible representation with highest weight $\mu$ and by
$\chi_\mu\in R(G)$ its character. Some additional notation 
to be introduced later: 
\begin{tabbing}
\indent \= $M_\mu$, $M_0=M\qu G$\quad \= \kill 
\> $J$ \> Weyl denominator; \ref{subsec:fix}\\
\> $\Alc$ \> fundamental alcove; \ref{subsec:fix}\\
\> $\rho$, $\alpha_0$, $c$ \> half-sum of positive roots, 
   highest root, dual Coxeter number; \ref{subsec:fix}\\
\> $B=B_k$ \> inner product on $\g$, 
   \ref{subsec:fix}\\
\> $\Lambda^*_k, R_k(G)$ \> level $k$ weights, level $k$ characters; 
   \ref{subsec:fix}\\
\> $T_k$ \> a certain finite subgroup of $T$; \ref{subsec:fix}\\
\> $t_\lam$ \> element of $T_{k+c}$ parametrized by $\lambda\in\Lambda^*_k$; 
   \ref{subsec:fix}\\
\> $G_\sig, LG_\sig$ \> 
   stabilizer in $G$ resp. $LG$ of face $\sig\subset \Alc$; \ref{sigplus}\\
\> $\mf{R}_{+,\sig}, \rho_\sig, \t_{+,\sig}$ \> 
   positive roots 
   for $G_\sig$, their half-sum, positive Weyl chamber
   ; \ref{sigplus}\\
\> $\gamma_\sig$ \> distinguished point in face $\sig\subset\Alc$; 
   \ref{subsec:loop}\\
\> $LG,\Om G$ \> free loop group, based loop group; \ref{subsec:loop}\\
\> $\Waff$ \> affine Weyl group; \ref{subsec:loop}\\
\> $\wh{LG}^{(k)}$ \> level $k$ central extension of the loop group; 
   \ref{subsec:central}\\
\> $\theta,\olt$    \> left, right Maurer-Cartan forms; \ref{subsec:loop}    \\
\end{tabbing}

\newpage 
\section{$\on{Spin_c}$-quantization of symplectic manifolds}
\label{sec:spinc}

In this section we review $\Spin_c$-quantization for compact
Hamiltonian $G$-manifolds. We explain how the fixed point formula for
the equivariant index, together with the ``quantization commutes with
reduction'' principle, leads to a formula \eqref{eq:formula} for the
index of a symplectic quotient in terms of fixed point contributions
for a certain finite subgroup.  This is the finite-dimensional version
of the main result of the paper.

\subsection{$\Spin_c$-quantization of Hamiltonian $G$-spaces}
We refer to Lawson-Michelson
\cite{la:sp} for  background on $\Spin_c$-structures, and to 
Duistermaat \cite{du:he} for a discussion of the symplectic
case.

Let $M$ be a compact, connected manifold with symplectic form $\om$,
together with a symplectic $G$-action. Given $\xi\in\g$ let $\xi_M=
\f{\p}{\p t}|_{t=0}\exp(-t\xi)^*$ denote the corresponding vector
field on $M$.  The action is called {Hamiltonian} if there exists a
$G$-equivariant map $\Phi\in C^\infty(M,\g^*)^G$ such that
$\iota(\xi_M)\om=\d\Phi(\xi)$ for all $\xi$. The map $\Phi$ is called
a moment map and the triple $(M,\om,\Phi)$ is called a Hamiltonian
$G$-space. By a theorem of Kirwan \cite{ki:con}, the intersection
$\Phi(M)\cap\t^*_+$ is a convex polytope; it is called the moment
polytope of $(M,\om,\Phi)$.

Suppose $(M,\om,\Phi)$ carries a $G$-equivariant pre-quantum line
bundle $L$. That is, $L$ comes equipped with an  
invariant connection $\nabla$, such that the Chern form
$c_1(L)=\f{i}{2\pi}\on{\curv}(\nabla)$ is equal to $\om$
and the vertical part of the vector field $\xi_L$ is given by 
Kostant's formula \cite{ko:qu}
\begin{equation}\label{Kostant}
 \on{Vert}(\xi_L)=2\pi\Phi(\xi)\, \f{\p}{\p \phi},
\end{equation}
where $\f{\p}{\p \phi}$ is the generator for the scalar $S^1$-action
on the fibers of $L$. 
Choose an invariant almost complex structure $I$ on $M$ that is
compatible with $\om$, in the sense that $\om(\cdot,I\cdot)$ defines a
Riemannian metric. 
The almost complex structure $I$ defines a
$G$-equivariant $\Spin_c$-structure on $M$, which we twist by the line
bundle $L$. 
Any choice of Hermitian connection on $TM$ defines a 
Dirac operator $\dirac$ for the twisted
$\Spin_c$-structure, and we define $\chi(M)\in R(G)$ as its equivariant
index
$$\chi(M)=\index_G(\dirac)\in R(G).$$
The index is independent of the choice of $I$ and  
of the connection. 
In case $M$ is K\"ahler and $L$ is holomorphic,  
the index coincides with the Euler characteristic for the sheaf
of holomorphic sections of $L$.

\subsection{Quantization commutes with reduction}
Following \cite{me:si} we call a point $\mu\in \g^*$ a {\em
quasi-regular} value of $\Phi$ if all $G_\mu$-orbits in
$\Phi^{-1}(\mu)$ have the same dimension. This includes regular values
and weakly regular values of $\Phi$.  For any
quasi-regular value $\mu\in\g^*$ the {\em reduced space} (symplectic
quotient) $ M_\mu=\Phinv(\mu)/G_\mu$
is a symplectic orbifold. If one drops the quasi-regularity
assumption, the space $M_\mu$ acquires more serious singularities
(cf. \cite{sj:st}).
For any dominant weight 
$\mu\in\Lambda^*_+$ which is a quasi-regular
value of $\Phi$,  
$$ L_\mu:=(L|_{\Phi^{-1}(\mu)}\otimes\C_{-\mu})/G_\mu\to M_\mu $$
is a pre-quantum orbifold-line bundle over $M_\mu$. The definition of 
$\Spin_c$-index carries over to the orbifold case, 
hence $\chi(M_\mu)$ is defined. In \cite{me:si}, this is extended
further to the case of singular symplectic quotients, using a partial
desingularization. The following Theorem was conjectured by
Guillemin-Sternberg and is known as ``quantization commutes with
reduction''.
\begin{theorem}[\cite{me:sym,me:si}] \label{qrc} Let $(M,\om,\Phi)$ be a
compact pre-quan\-tized Ha\-mil\-to\-ni\-an $G$-ma\-ni\-fold.  Then the
multiplicity of $\chi_\mu$ in $\chi(M)\in R(G)$ is equal to
$\chi(M_\mu)$.
\end{theorem} 
In particular, only weights $\mu\in\Lambda^*_+$ which
are contained in the moment polytope $\Phi(M)\cap\t^*_+$ can appear in
$\chi(M)$. 

\subsection{Equivariant index theorem}
\label{sec:fixed}
The equivariant index theorem expresses the value $\chi(M,g)$
in terms of local data at the fixed point set $M^g$. 
We recall that the
connected components $F\subseteq M^g$ are compact, 
embedded almost complex 
submanifolds of $M$. They are invariant under the action of the centralizer
$G_g$, and the pull-backs $\om_F,\Phi_F,L_F$ 
of $\om,\Phi,L$ give $F$ the structure of a pre-quantized Hamiltonian 
$G_g$-space. The action of
$g\in G$ on $L$ restricts to a multiplication by a phase factor
$\mu_F(g)\in\U(1)$ on $L_F$. 

Let $\Td(F)$ be the Todd form, for any
Hermitian connection on $TF$, and let the form 
$\D_\C(\nu_F,g)$ be defined by 
$$ \D_\C(\nu_F,g)={\det}_\C(1-A_F(g)^{-1}
e^{R_{\nu_F}/2\pi}).
$$ 
Here $A_F(g)\in\Gamma^\infty(\U(\nu_F))$ is
the unitary bundle automorphism of $\nu_F$ induced by $g$, and 
$R_{\nu_F}\in\Om^2(F,\mf{u}(\nu_F))$ the 
curvature of an invariant Hermitian
connection on $\nu_F$. 
The Atiyah-Segal-Singer fixed point formula 
\cite{at:1,at:2} 
asserts that
\begin{equation}
\label{fp1}
\chi(M,g)=\sum_{F\subseteq M^g}\chi(\nu_F,g)
\end{equation}
where 
\begin{equation}
\label{fp2}
\chi(\nu_F,g)=\mu_F(g)\,\int_F\f{\Td(F)e^{c_1(L|_F)}}{\D_\C(\nu_F,g)}.
\end{equation}
(For finite fixed point set, the formula 
is a special case of the Atiyah-Bott Lefschetz formula 
\cite{at:le2}.) We will also need an equivalent expression for 
the fixed point contributions, in which the almost complex structure 
enters only via the $\Spin_c$-line bundle $\L$, given as a tensor product 
$\L=L^2\otimes K^{-1}$ of $L\otimes L$ with the 
anti-canonical line bundle $K^{-1}$ for the almost 
complex structure:
\label{sec:alter}
\begin{equation}\label{lco}
 \chi(\nu_F,g)=\zeta_F(g)^{1/2}\, \int_F \f{\hat{A}(F) e^{\hh
c_1(\ca{L}|_F)}}{\D_\R(\nu_F,g)}.
\end{equation}
Here $\zeta_F(g)$ is the eigenvalue for the action of $g$ on 
$\L|_F$, and the square root $\zeta_F(g)^{1/2}=\mu_F(g)\kappa_F(g)^{-1/2}$
is defined as 
$$ \kappa_F(g)^{-1/2}:=\det(A_F(g)^{1/2}).$$
where $A_F(g)^{1/2}\in \Gamma^\infty(U(\nu_F))$ 
is the unique square root of $A_F(g)$ having all its eigenvalues in the set 
$\{e^{i\phi}|\,0\le \phi<\pi\}$. 
Furthermore $\wh{A}(F)$ is the $\hat{A}$-form of $F$, and 
$$  
\D_\R(\nu_F,g)= i^{\on{rank}(\nu_F)/2}
{\det}^{1/2}_\R(1-A_F(g)^{-1}e^{R_{\nu_F}/2\pi}),
$$
viewing $A_F(g)$ as a real automorphism, 
$R_{\nu_F}$ as a $\mf{o}(\nu_F)$-valued 2-form, and 
taking the positive square root. 
The fixed point expressions \eqref{fp2} and \eqref{lco}
are identical because 
$$ \D_\C(\nu_F,g)=\D_\R(\nu_F,g)e^{\hh c_1(K_{\nu_F})} \kappa_F(g)^{1/2}$$
and 
$$ \Td(F)=\hat{A}(F)e^{-\hh c_1(K_F)}$$
where $K_F$ is the canonical bundle of $F$ and $K_{\nu_F}$ that 
for $\nu_F$.

\subsection{Fixed point formula for multiplicities}
\label{subsec:fix}

We now use finite Fourier transform to extract the multiplicity of
any weight $\mu\in\Lambda^*_+$ from Formula \eqref{fp1}.  
For a different approach using Fourier series, see
Guillemin-Prato \cite{gu:he}. For simplicity, we only 
consider the case where $G$ is simply connected.  

We need to introduce some extra 
notation and facts regarding compact Lie groups.  
Suppose (for a short moment) that $G$ is simple.
Let $\alpha_0\in\Lambda^*$ be the highest root, $h_{\alpha_0}\in\t$
its coroot, $\rho\in\Lambda^*$ the half-sum of positive roots, and
$c=1+\rho(h_{\alpha_0})$ the dual Coxeter number. The 
fundamental alcove is denoted
$\Alc:=\{\xi\in \t_+|\alpha_0(\xi)\le 1\}$. Let the {\em basic 
inner product} $B^G$ be the unique invariant inner product on $\g$ such 
that $B^G(h_{\alpha_0},h_{\alpha_0})=2$. It has the important property 
that it restricts to an integer-valued $\Z$-bilinear form on 
the lattice $\Lambda$ (see \cite{br:rep}, Chapter V.2). 

For a general simply
connected group $G$, decompose into simple factors $G=G_1\times
\ldots\times G_s$ with dual Coxeter numbers $c=(c_1,\ldots,c_s)$, and
define $\Alc=\Alc_1\times\ldots\times \Alc_s$. Any invariant 
symmetric bilinear form on $\g$ can be written 
$$ B_k:= \sum_{j=1}^s k_j B^{G_j}, $$
where $k_j\in\R$. We denote by $B_k^\flat:\,\g\to \g^*$ the linear
map defined by $B_k$, and if all $k_j\not=0$ the inverse map is 
denoted $B_k^\sharp=(B_k^\flat)^{-1}$. Suppose 
all $k_j$ are positive integers. Then $B_k$ is integer-valued
on $\Lambda$, and we have inclusions 
$B_k^\flat(\Lambda)\subset\Lambda^*$ and 
$B_k^\sharp(\Lambda^*)\supset\Lambda$.
The finite Fourier transform is taken using the finite subgroup
$$T_k:=B_k^\sharp(\Lambda^*)/\Lambda$$ 
of $T=\t/\Lambda$. 
Let 
$$ \Alc^*_k:=B_k^\flat(\Alc)\subset\t^*_+,\ \
\Lambda^*_k=\Lambda^*\cap\Alc^*_k.$$
The weights in $\Lambda^*_k$ are
called {\em weights at level $k$}. Using the definition of the alcove, one
verifies that the map
$$ \Lambda^*_k\to T_{k+c},\
\lambda\mapsto t_\lambda:=\exp(B_{k+c}^\sharp(\lambda+\rho)) $$
takes values in $T_{k+c}^{\on{reg}}=T_{k+c}\cap G^{\on{reg}}$ and identifies 
$\Lambda^*_k=T^{\on{reg}} _{k+c}/W$. 
\begin{example}
Suppose $G=\SU(2)$. Then $\rho\in\t^*$ spans the weight lattice, and
$\alpha_0=2\rho$ is the positive root. Therefore,
$B_k^\flat(\Lambda)=2k\Lambda^*$ and $T_k=\Z_{2k}$.  Furthermore,
$\Alc^*_k=\{t\rho|\,0\le t\le k\}$ and
$\Lambda^*_k=\{0,\rho,\ldots,k\rho\}$. The dual Coxeter number is
$c=2$. For $\lambda=l\rho$, the element $t_\lambda$ is a diagonal
matrix with entries $e^{i\phi}, e^{-i\phi}$ down the diagonal, where
$\phi=\pi \f{l+1}{k+2}$.
\end{example}
Let the set of level $k$ characters $R_k(G)\subset R(G)$ 
be the additive subgroup generated by all $\chi_\mu$ with 
$\mu \in \Lambda^*_k$. 
One has the orthogonality relations,
\begin{eqnarray}
\label{eq:Mul}
\sum_{\lambda\in\Lambda^*_k} |\JJ (t_\lambda)|^2\,
\chi_\mu(t_\lambda)\chi_{\mu'}(t_\lambda)^*&=&
{\# T_{k+c}}\,\,
\delta_{\mu,\mu'},\ \ \ 
\mu,\mu'\in\Lambda^*_k\\ 
\sum_{\mu \in\Lambda^*_k} |\JJ (t_\lam)|^{2}
\chi_\mu(t_\lam)\chi_\mu(t_{\lam'})^*&=&
{{\# T_{k+c}}}\,\,
\delta_{\lam,\lam'},\ \ \ 
\lam,\lam'\in \Lambda^*_k
\label{eq:Mul2}
\end{eqnarray}
where $\JJ :\,T\to \C$ is the Weyl denominator
$$\JJ (t)=\sum_{w \in W} (-1)^{l(w)} t^{w\rho}.$$
%
%
%
These formulas are obtained from the Weyl character formula and finite
Fourier transform for $T_{k+c}\subset T$. As a consequence, 
every level $k$ character is determined by its restriction to 
$T_{k+c}^{reg}$. One may view $R_k(G)$ as a quotient of $R(G)$ 
by the ideal of characters vanishing on $T_{k+c}^{reg}$; 
this defines a ring structure on $R_k(G)$ known as {\em fusion product}.

Let us return to the problem of calculating multiplicities in
$\chi(M)$ from the fixed point formula.  
The remark following Theorem \ref{qrc} shows that 
$\chi(M)\in R_k(G)$ provided $k_j>0$ are chosen large enough  
so that $\Alc^*_k$ contains the moment polytope $\Phi(M)\cap\t^*_+$. 
The multiplicity of any
$\mu\in\Lambda^*_k$ can be computed from the orthogonality relation
\eqref{eq:Mul}, substituting the Atiyah-Segal-Singer fixed point
formula for $\chi(M,t_\lambda)$. On the other hand this multiplicity
equals $\chi(M_\mu)$, by Theorem \ref{qrc}. This gives,
\begin{proposition}  Let $(M,\om,\Phi)$ be a pre-quantized 
compact Hamiltonian $G$-manifold with character $\chi(M)$ at level $k$.
Then the index of the reduced space $M_\mu$ can be expressed in terms of fixed
point contributions in $M$:
\begin{equation}\label{eq:formula}
\chi(M_\mu)=  \frac{1}{{\# T_{k+c}}} \sum_{\lam \in \Lambda^*_k}
\chi_{\mu}(t_\lam)^*
|\JJ (t_\lam)|^2 
\sum_{\ F\subseteq M^{t_\lambda}} \chi(\nu_F,t_\lam) .
\end{equation}
\end{proposition}

\section{Review of loop group actions and group-valued moment
maps}
\label{sec:review}

\subsection{Loop groups}
\label{subsec:loop}
For the material in this subsection we refer to Pressley-Segal
\cite[Section 4.3]{pr:lo}. Let $S^1=\R/\Z$ be the parametrized circle
with coordinate $s$. 

Throughout we will fix a ``Sobolev level'' $f>1$. 
We denote by $\Oma^0(S^1,\g)$ the space of $\g$-valued
$0$-forms of Sobolev class $f+1/2$ and by 
$\Oma^1(S^1,\g)$ the $\g$-valued $1$-forms of Sobolev class $f-1/2$.
Then forms in $\Oma^0(S^1,\g)$ are $C^1$ and those 
in $\Oma^1(S^1,\g)$ are $C^0$. 
Let the {\em (free) loop group} $LG$ consist of maps 
$S^1 \to G$ of Sobolev class $f+1/2$. It is a Banach Lie group 
with Lie algebra $L\g={\Oma}^0(S^1,\g)$. The kernel of the
evaluation mapping $LG\to G,\,g\mapsto g(0)$ is called the {\em based loop
group} $\Om G$.  The free loop group  is a semi-direct product 
$LG=G\ltimes \Om G$ where $G$ is embedded as constant loops and the
action of $G$ on $\Om G$ is pointwise conjugation.  We embed the
lattice $\Lambda\subset\t$ into $LG$ by the map which takes
$\xi\in\Lambda$ to the loop,
$$
\R/\Z\to G,\ \ s\mapsto \exp(-s\xi).
$$
%

The loop group $LG$ acts on the affine space $\A(S^1)=\Oma^1(S^1,\g)$
of connections on the trivial $G$-bundle over $S^1$  
by gauge transformations 
$$
g \cdot \mu = \Ad_g(\mu)-g^*\olt.
$$
Here $\ol{\theta} \in \Om^1(G,\g)$ denotes the right-invariant
Maurer-Cartan form.
The orbit space for this action can be described as follows. 
Consider the embedding $\t\subset\g\hra {\Oma}^1(S^1,\g)$ 
by the map $\xi\mapsto \xi\d s$. 
The intersection of any  
$LG$-orbit with $\t$ is an 
orbit of the affine Weyl group $\Waff=W\ltimes \Lambda$. 
Hence there are natural identifications, 
$${\Oma}^1(S^1,\g)/LG=\t/\Waff=T/W=G/\Ad(G).$$
In particular, one has a 1-1 correspondence between $LG$-orbits and
conjugacy classes in $G$. For $G$ simply connected, all of these sets
are also identified with the fundamental alcove $\Alc$. That is, each
coadjoint $LG$-orbit meets the alcove $\Alc$ in exactly one point.

For $s\in\R$ and any $\mu\in\Oma^1(S^1,\g)$ let 
$\Hol_s(\mu)\in G$ denote the parallel transport from $0$ to $s$. 
Thus $s\mapsto h(s)=\Hol_s(\mu)$ is the unique solution of 
the initial value problem
$h'(s)h(s)^{-1}\equiv h^*\olt=\mu,\ h(0)=e$. One has the 
equivariance property 
\begin{equation}\label{eq:hols}
\Hol_s(g\cdot\mu)=g(s)\Hol_s(\mu)g(0)^{-1},
\end{equation} 
showing in particular that the based gauge group $\Om G$ acts {\em
freely}.  
Let $\Hol(\mu):= \Hol_1(\mu)$ denote holonomy of $\mu$
around $S^1$. Any two elements $\mu_1,\mu_2 \in \Oma^1(S^1,\g)$ with
$\Hol(\mu_1) = \Hol(\mu_2)$ are related by a based gauge
transformation, $\Hol_s(\mu_1)\Hol_s(\mu_2)^{-1}$.  The {\em holonomy
map}
$$ \Hol:\,\Oma^1(S^1,\g)\to G$$
gives $\Oma^1(S^1,\g)$ the structure of a Banach principal $\Om
G$-bundle over $G$.  It is equivariant with respect to the evaluation
map $LG\to G$ and once again gives the correspondence between
$LG$-orbits and conjugacy classes.

For any $\mu\in {\Oma}^1(S^1,\g)$, the evaluation map $LG\to G$ 
induces an isomorphism, $LG_\mu\cong  G_{\Hol(\mu)}$, with inverse map
\begin{equation}\label{eq:inverse}
 G_{\Hol\mu}\to LG_\mu,\ \ g\mapsto \Ad_{\Hol_s(\mu)}(g).
\end{equation}
In particular this shows that all stabilizer groups $LG_\mu$ are 
compact. 

\label{sigplus}
If we make the additional assumption that $G$ is simply connected,
then all centralizers $G_g$, hence also all stabilizer groups
$LG_\mu\cong G_{\on{Hol}(\mu)}$, are connected.  For any open face
$\sig\subset\Alc$, the stabilizer group $LG_\mu$ of $\mu\in\sig$ is
independent of $\mu$ and will therefore be denoted $LG_\sig$.  The
evaluation map defines an isomorphism $LG_\sig\cong G_\sig$ with the
centralizer of $g=\exp(\mu)$. If $\sig\subset\ol{\tau}$ then
$LG_\tau\subset LG_\sig$ and $G_\tau\subset G_\sig$; in particular
every $G_\sig,\,LG_\sig$ contains the maximal torus $T$.  The root
system $\mf{R}_{\sig}$ of $G_\sig$ consists of all $\alpha\in\mf{R}$
such that the restriction $ \alpha |_\sig$ is integer valued.  From
the definition of the dual Coxeter number, it follows that
$B_c^\sharp(\rho)\in\on{int}(\Alc)$. Let $\mf{R}_{+,\sig}$ be all
$\alpha\in\mf{R}_\sig$ such that $\alpha(B_c^\sharp(\rho)) \ge
\alpha |_\sig$, and let $2\rho_\sig$ be their sum.  One can check
(cf. \cite{me:can}) that
\begin{equation} \label{eq:gammasig}
\gamma_\sig:=B_c^\sharp(\rho-\rho_\sig)
\end{equation} 
is always contained in $\sig$. The positive Weyl chamber 
$\t_{+,\sig}$ for $\mf{R}_{+,\sig}$ is the cone over 
$\Alc-\gamma_\sig$.\label{rhosig}

\subsection{Hamiltonian loop group actions}
The space $L\g^*:=\Oma^1(S^1,\g^*)$ is a dense sub-space of the
topological dual space of $L\g={\Oma}^0(S^1,\g)$, using the pairing of
$\g^*$ and $\g$ followed by integration over $S^1$.  Given an
invariant inner product $B$ on $\g$, the isomorphism
$B^\flat:\,\g\to\g^*$ gives rise to an identification $B^\flat:\,
\Oma^1(S^1,\g)\to L\g^*$. The affine $LG$-action induced on $L\g^*$
via this isomorphism will be called the {\em coadjoint loop group
action}, and its orbits will be called coadjoint $LG$-orbits.  
Recall that a 2-form on a Banach
manifold $\wh{M}$ is called weakly symplectic if its kernel is trivial
everywhere.
\begin{definition}  
A Hamiltonian $LG$-manifold is a triple $(\widehat{M}, \widehat{\om},
\widehat{\Phi})$, consisting of a Banach manifold
$\widehat{M}$ with a smooth $LG$-action, 
a weakly symplectic 2-form $\wh{\om}$
on $\wh{M}$, and an equivariant map $\widehat{\Phi}: \ \wh{M}
\to L\g^*$ satisfying 
\begin{equation}\label{eq:lpmoment}
\iota(\xi_{\wh{M}})\wh{\om} 
= \d \widehat{\Phi}(\xi),\ \ \xi\in L\g .
\end{equation}
\end{definition}
Much of the theory of compact Hamiltonian $G$-spaces carries over 
to Hamiltonian loop group spaces if one assumes that the moment map
$\wh{\Phi}$ is {\em proper}. For example, if $\mu\in L\g^*$ is a 
regular value of the moment map then the reduced space 
$\wh{M}_\mu=\wh{\Phi}^{-1}(\mu)/LG_\mu$ is a compact, 
finite dimensional symplectic orbifold. More generally, this holds
true for quasi-regular values $\mu$, i.e. if all $LG_\mu$-orbits in
$\wh{\Phi}^{-1}(\mu)$ have the same dimension. 

For $G$ simply connected and $\wh{M}$ connected, 
there is also a convexity theorem 
\cite[Theorem 4.11]{me:lo}, stating the intersection
$\wh{\Phi}(\wh{M})\cap\Alc$ is a convex polytope. 
We refer to this polytope as the moment polytope of $\wh{M}$.

Basic examples for Hamiltonian $LG$-spaces with proper moment maps
are coadjoint orbits $\O=LG\cdot\mu$ for $\mu\in L\g^*$, with 
moment map the inclusion. The 2-form is uniquely determined
by the moment map condition, and is given by an analog to the  
Kirillov-Kostant-Souriau formula. 
The motivating examples of moduli spaces of flat connections on
surfaces with boundary are discussed in the companion paper
\cite{al:ve}.

\subsection{Group-valued moment maps}
Suppose $(\wh{M},\wh{\om},\wh{\Phi})$ is a Hamiltonian $LG$-space with
proper moment map. As mentioned above, the holonomy map
$\Hol:\,L\g^*\to G$ is a principal $\Om G$-bundle.  By equivariance of
the moment map, the based loop group $\Om G$ also acts freely on
$\wh{M}$, and by properness the quotient $M:=\wh{M}/\Om G$ is a
compact, smooth, finite dimensional manifold \cite[Section
3.2.1]{me:co}. The action of $LG=G\ltimes \Om G$ descends to a
$G$-action on $M$, and the moment map $\wh{\Phi}$ to a $G$-equivariant
map $\Phi:\,M\to G$, which makes the following diagram commute:
$$ \begin{array}{ccc} \widehat{M} & \stackrel{\widehat{\Phi}}{\lra} &
L\g^* \\ \downarrow & & \downarrow \\ M & \stackrel{\Phi}{\lra} & G
\end{array}
$$
We call $M$ the holonomy manifold of $\wh{M}$.  In \cite{al:mom} $M$
is interpreted as a Hamiltonian $G$-space with group valued moment map
$\Phi$. The definition is as follows. Let $\theta,\ol{\theta}$ be the
left resp. right invariant Maurer-Cartan forms on $G$, and 
$\eta\in\Om^3(G)$ the
canonical closed, bi-invariant 3-form
$$ \eta = \frac{1}{12} B(\theta, [\theta,\theta])=
\frac{1}{12} B(\olt, [\olt,\olt]).$$
Define a 2-form on $L\g^*$ by 
$$ 
\varpi=\hh \int_0^1 B( \Hol_s^*\olt, \f{\p}{\p s}\Hol_s^*\olt )\ 
\d s.$$
By \cite[Proposition 8.1]{al:mom} the form $\varpi$ has the 
property  $ \d\varpi=\Hol^*\eta$, and its contractions with 
generating vector fields for the 
$LG$-action are 
$$ \iota(\xi_{L\g^*})\varpi = -\d \mu(\xi)
+ \hh \Hol^*\,B(\theta+\olt,\xi(0)). 
$$
Here $\mu(\xi)$ is the function on $L\g^*$ taking $\mu\in L\g^*$ 
to the pairing with $\xi\in L\g$. It follows that 
$\wh{\om}-\wh{\Phi}^*\varpi\in\Om^2(\wh{M})=\Hol^*\om$ 
for a unique 2-form $\om\in\Om^2(M)$. As shown in 
\cite{al:mom}, this 2-form has properties 
\begin{equation}\label{eq:1}
\d\om=\Phi^*\eta,
\end{equation}
\begin{equation}\label{eq:2}
\iota(\xi_M)\om=\hh \Phi^*B(\theta+\olt,\xi),\ \ \ \ 
\mbox{ for all }\xi\in\g,
\end{equation}
\begin{equation}\label{eq:3}
\ker(\om_m)=\{\xi_M(m)|\ \Ad_{\Phi(m)}\xi=-\xi\}
\end{equation}
and conversely every compact $G$-manifold $M$ with an 
equivariant map $\Phi\in C^\infty(M,G)$ and an invariant 2-form 
$\om$ with these three properties defines a Hamiltonian 
$LG$-manifold with proper moment map. We call $(M,\om,\Phi)$ with properties 
\eqref{eq:1}, \eqref{eq:2}, \eqref{eq:3} a {\em group-valued 
Hamiltonian $G$-space}.  
As for $\g^*$-valued moment maps, we call an element 
$g\in G$ a quasi-regular value of $\Phi$ if all 
$G_g$-orbits in $\Phi^{-1}(g)$  have the same dimension. 
It is proved in \cite{al:mom} that in this case the reduced 
space $M_g=\Phi^{-1}(g)/G_g$ is a symplectic orbifold, 
with 2-form induced from $\om$. Given $\mu\in L\g^*$ 
with $g=\Hol(\mu)$, one finds that $g$ is quasi-regular
for $\Phi$ if and only if $\mu$ is quasi-regular for $\wh{\Phi}$, 
and $\wh{M}_\mu=M_g$ as symplectic spaces.

\section{The fixed point formula for loop group actions}
\label{sec:mainresult}
In this Section, we discuss pre-quantization of Hamiltonian loop
group manifolds and state the main result of this paper, 
the fixed point formula Theorem \ref{onal}.

\subsection{Central extension of $LG$}
\label{subsec:central}
We will assume for the rest of this paper that $G$ is simply connected, 
and let $G=G_1\times \ldots \times G_s$ be its decomposition into simple factors.  
Each of the bilinear forms $B=B_l$, $l\in \R^s$, on $\g$ defines a central
extension $\wh{L\g}$ of $L\g=\Oma^0(S^1,\g)$ by $\R$, with cocycle
$$\mf{c}:\,L\g\times L\g \to \R,\ 
\mf{c}(\xi_1,\xi_2)=\int_{S^1}\,B(\xi_1,\d\xi_2).$$ 
The dual action of $LG$ on $\wh{L\g}^*={\Oma}^1(S^1,\g^*)\times\R$
is given by the formula, $g\cdot (\mu,\tau)=(\mu-\tau\,
B^\flat(g^*\olt),\tau)$. If $B$ is non-degenerate, 
it identifies $\wh{L\g}^*={\Oma}^1(S^1,\g)\times\R$ 
and the gauge action of $LG$ on $\A(S^1)$ becomes the 
action on the affine hyperplane $\tau=1$. 

It is known \cite[Theorem (4.4.1)]{pr:lo} that the Lie algebra 
extension exponentiates to a group extension 
\begin{equation} \label{eq:grex}
1\to \U(1)\to \wh{LG}\to LG\to 1
\end{equation}
exactly if all $l_j$ are integers. Since $LG$ is connected and simply 
connected, the extension is unique. For any subgroup $\H \subset LG$, 
we denote by $\wh{\H}$ the pull-back of the central extension. Since the 
defining cocycle $\mf{c}$ vanishes on $\g$, and since $G$ is connected and 
simply connected, the central extension 
$\wh{G}$ is canonically trivial.  That is, there is an embedding
\begin{equation} \label{eq:triv}
G \to \wh{LG},
\end{equation}
and $\wh{LG}$ is a semi-direct product $\wh{LG}=G\ltimes \wh{\Om G}$. 
The following proposition describes the central extensions for 
various subgroups of $LG$.  
\begin{proposition}\label{prop:resext}
Let $l\in \Z^s$ and $\wh{LG}=\wh{LG}^{(l)}$ the corresponding 
central extension of the loop group. 
\begin{enumerate}
\item 
Let $\xi,\xi'\in \Lambda$, and $\hat{\xi},\hat{\xi}'\in\wh{\Lambda}$ 
arbitrary lifts. Then the Lie group commutator is given by the formula
$$ [\hat{\xi},\hat{\xi}']=(-1)^{B_l(\xi,\xi')}.$$
In particular, if all $l_j$ are even, the central extension
$\wh{\Lambda}$ is trivial. 
\item
Embed $T \hra \wh{LG}$ using \eqref{eq:triv}. The central extension of 
$LT\subset LG$ is a semi-direct product 
$$\wh{LT}=T\ltimes\wh{\Om T}$$
where the action of $T$ on $\wh{\Om T}$, 
over the connected component of 
$\Om T$ containing $\xi\in\Lam$, is given by multiplication 
with $t^{B^\flat(\xi)}$.
\item
Suppose all $l_j>0$.
The central extension of $T_l\times \Om T$ is a direct product, 
$$ \wh{T_l\times \Om T}=T_l\times\wh{\Om T}.$$ 
\end{enumerate}
\end{proposition}

\begin{proof}
Part (a) is proved by Pressley-Segal \cite[Section 4.8]{pr:lo} 
for simply laced groups, and by Toledano Laredo 
\cite[Proposition 3.1]{to:po} in the general case. 
Part (c) follows from part (b) since 
$t^{B^\flat(\xi)}=1$ for $\xi\in \Lambda$, $t\in T_l$. 

It remains to prove part (b).  Since the conjugation action of $T$ on
$\Om T$ is trivial, its action on any connected component of 
$\wh{\Om T}$ is scalar multiplication on the fibers 
by some character for $T$. To compute the weight for this action at 
$\xi\in \Lambda\subset \Om T$, let $\alpha\in \Om^1(\wh{LG})$ be
the left-invariant connection 1-form defined by the splitting
$\wh{L\g}=L\g\times\R$. The weight $\mu\in\Lambda^*$ for the
$T$-action on the fiber over $\xi$ is given by
$$\l\mu,\zeta\r=\l \alpha_{\hat{\xi}},\zeta_{\wh{LG}}\r,\ \ \zeta\in\t$$
where $\hat{\xi}$ is any lift of $\xi$. In left trivialization of 
the cotangent bundle of $\wh{LG}$, 
$\alpha$ is the constant map from $\wh{LG}$ to 
$(0,1)\in\wh{L\g}^*$. This shows that its contraction 
with the left-invariant vector field generated by 
$\zeta$ is zero. To compute its contractions with the right-invariant 
vector field generated by $\zeta$, we note that under the right action of $\hat{\xi}$,   
$$R_{\hat{\xi}^{-1}}^*\alpha=\Ad_{\hat{\xi}}^*\alpha=
\hat{\xi}\cdot (0,1)=(B^\flat(\xi),1).$$
Since $\zeta_{\wh{LG}}$ is the difference 
between left and right invariant vector fields, we find that 
$\l \alpha_{\hat{\xi}},\zeta_{\wh{LG}}\r=
\l B^\flat(\xi),\zeta\r$, proving $\mu=B^\flat(\xi)$.
\end{proof}
Below we will often use the following terminology.
A {\em level $l$ line bundle} over an $LG$-space $X$ is 
an $\wh{LG}^{(l)}$-equivariant 
line bundle $L\to X$ where the central $\U(1)$ acts by 
scalar multiplication with weight 1. Equivalently, $L$
carries actions of the central extensions $\wh{LG_j}^{(1)}$
where the central $\U(1)$'s act with weights $l_j$. The tensor product
of two line bundles at levels $l,l'$ is at level $l+l'$, and line
bundles at level $0$ are simply $LG$-equivariant line bundles.
%

\subsection{Fixed point formula for Hamiltonian loop group actions}
Suppose now that $B=B_k$ where all $k_j$ are positive integers, and
let $\wh{LG}=\wh{LG}^{(k)}$.  Using $B$ we identify $\g\cong\g^*$ and
$\Oma^1(S^1,\g)\cong L\g^*$.  The formula for the coadjoint action of
$LG$ on $\wh{L\g}^*$ shows that a Hamiltonian $LG$-manifold
$(\wh{M},\wh{\om},\wh{\Phi})$ is equivalently a Hamiltonian
$\wh{LG}$-manifold on which the central $\U(1)$ acts trivially, with
moment map $1$.  Accordingly we define a pre-quantum line bundle to be
a level $k$ line bundle $\nh{L}\to\wh{M}$ with invariant connection, 
such that the Chern form $c_1(\nh{L})$ is equal to the
symplectic form $\wh{\om}$ and such that the vertical part of the
fundamental vector fields is given by Kostant's formula.
\begin{example}
The coadjoint orbit $LG\cdot\mu$ through $\mu\in \Alc^*_k$ is 
pre-quantizable if and only if 
$\mu\in\Lambda^*_k$. The pre-quantum line bundle is the 
associated bundle $\wh{LG}\times_{\wh{LG}_\mu}\C_{(\mu,1)}$. 
If $\Sig$ is a compact oriented surface with 
boundary, the moduli space $\M(\Sig)$ of flat $G$-connections 
by based gauge equivalence is pre-quantizable. For non-simply 
connected groups the situation is more complicated -- see the 
companion paper \cite{al:ve}.
\end{example}
Suppose the moment map $\wh{\Phi}$ is proper and 
$\mu\in\Lambda^*_k$ is a quasi-regular value. Then 
$\wh{M}_\mu$ is a symplectic orbifold with pre-quantum 
orbifold line bundle 
$$
\nh{L}_\mu=(\nh{L}|_{\wh{\Phi}^{-1}(\mu)}\otimes\C_{-(\mu,1)})/
\wh{LG}_\mu.
$$
(Notice that in fact $LG_\mu$ acts on the tensor product, for the
central $\U(1)\subset \wh{LG}_\mu$ acts with weights $+1$ and $-1$ on
the two factors and hence acts trivially on the product.) Hence the
$\Spin_c$-index $\chi(\wh{M}_\mu)$ is defined for quasi-regular
$\mu\in\Lambda^*_k$. If $\mu$ is not quasi-regular one can still
define $\chi(\wh{M}_\mu)$ by a partial desingularization as in
\cite{me:lo}.
The main result of the paper is the following loop group analog 
of Equation \eqref{eq:formula}.  
\begin{theorem}[Fixed Point Formula] \label{onal} 
Suppose $(\wh{M},\wh{\om},\wh{\Phi})$ is a pre-quantized 
Hamiltonian $LG$-manifold with proper moment map. 
For all level $k$ weights $\mu\in\Lambda^*_k$, the index of 
the symplectic quotient $\wh{M}_\mu$ is a sum of integrals 
over fixed point manifolds $F\subseteq M^t$ for the 
action of elements $t\in T_{k+c}^{reg}$ on the holonomy 
manifold $M=\wh{M}/\Om G$,
$$ \chi(\wh{M}_\mu) = \frac{1}{{\# T_{k+c}}} \sum_{\lam \in
\Lambda^*_k}
\chi_\mu(t_\lam)^* |\JJ (t_\lam)|^2   \sum_{F \subseteq M^{t_\lam}}
\zeta_F(t_\lam)^{1/2}
\int_F \frac{\widehat{A}(F) e^{\hh c_1(\L_F) }}{ \D_\R(\nu_F,t_\lam)} . 
$$
The terms entering the fixed point contributions will be explained
in Subsection \ref{sec:fixedpointcontributions} below.
\end{theorem}
Equivalently, organizing the $\Spin_c$-indices 
$\chi(\wh{M}_\mu)$ into a level $k$ character 
$\chi(M):=\sum_{\mu\in\Lambda^*_k}\chi(\wh{M}_\mu) \chi_\mu$
we have the formula 
\begin{equation} \label{onaleq}
\chi(M,t) = \sum_{F\subseteq M^t}
\chi(\nu_F,t),\ \ \ \ t\in T_{k+c}^{reg}
\end{equation}
where 
\begin{equation} \label{eq:LF}
 \chi(\nu_F,t) = \zeta_F(t)^{1/2}
\int_F \frac{\widehat{A}(F) e^{\hh c_1(\L_F) }}{ \D_\R(\nu_F,t)} . 
\end{equation}
The proof of Theorem \ref{onal} will be given in the final 
Section of this paper. 
%
\subsection{The fixed point contributions}
\label{sec:fixedpointcontributions}
In general, the holonomy manifold $M$ does not carry a naturally 
induced $\Spin_c$-structure, even though 
the expressions \eqref{eq:LF} resemble the fixed point contributions
of a $\Spin_c$-Dirac operator. Our strategy for defining the terms
entering \eqref{eq:LF} is to first restrict data on $\wh{M}$ to 
a certain finite dimensional submanifold $\wt{F}\subset   \wh{M}^t$
covering $F$, and then to show that the restrictions descend to $F$
itself.
\begin{proposition}\label{prop:geometry}
Let $(M,\om,\Phi)$ be a Hamiltonian $G$-space with group valued moment
map, and $(\wh{M},\wh{\om},\wh{\Phi})$ the corresponding loop group space. 
Let $t\in T^{\on{reg}}$, and $F\subseteq M^t$ a connected component of 
the fixed point set. Let $\wh{F}$ be the pre-image of $F$ under the map 
$\wh{M}^t\to M^t$, and $\wt{F}$ the intersection of $\wh{F}$ 
with $\wh{\Phi}^{-1}(\t)$.  
\begin{enumerate}
\item
$F$ is a group-valued Hamiltonian $T$-space, with symplectic form 
$\om_F$ and moment map $\Phi_F$ the pull-backs of ${\om},{\Phi}$.
\item
$\wh{F}$ is a (possibly disconnected) Hamiltonian $LT$-manifold, with 
2-form and moment map the pull-backs of $\wh{\om},\wh{\Phi}$. 
It has $(F,\om_F,\Phi_F)$ as its holonomy manifold.
\item
$\widetilde{F}$ is a finite-dimensional Hamiltonian $T$-manifold, with 
2-form and moment map the pull-backs of $\wh{\om},\wh{\Phi}$. It 
carries a free symplectic action of the lattice 
$\Lambda$ which commutes with the 
action of $T$. One has $F=\widetilde{F}/\Lambda$ as a symplectic 
$T$-manifold, and $\widetilde{F}=\wh{F}/\Om_0 T$ where $\Om_0 T$ is 
the identity component of $\Om T$. 
\end{enumerate}
\end{proposition}

\begin{proof}
We begin by showing that $F$ is symplectic.  Since $t$ is a regular
element, $\Phi(F)\subseteq G^t=T$. 
Let $m\in F$, $g=\Phi(m)$, and consider the splitting 
of the tangent space
$$ T_mM=E\oplus \g_g^\perp$$ 
where $E=(\d_m\Phi)^{-1}(\g_g)$ and the second summand is embedded by
the generating vector fields. By \cite[Section 7]{al:mom}, 
the splitting is $\om$-orthogonal and the 
restriction of $\om$ to $E$ is symplectic. Since the 
action of $t$ on $E$ preserves 
the 2-form, the subspace $T_mF=(T_mM)^t=E^t$ is symplectic as well.
This shows that $\om_F$ is non-degenerate. It is closed since
$\d\om_F=\iota_F^*\Phi^*\eta=\Phi_F^*\iota_T^*\eta=0$. 
The moment map condition for $(F,\om_F,\Phi_F)$ follows from that for
$(M,\om,\Phi)$. This proves (a). 

Clearly $\wh{F}$ is $\Om T$-invariant, and its 
image under the moment map is contained in 
$(L\g^*)^t=L\t^*$. Also $\wh{F}/\Om T\subseteq F$. 
To prove the reverse inclusion, suppose $m\in F$. 
Its pre-image under the map $\wh{M}\to M$ meets 
$\wh{\Phi}^{-1}(L\t^*)$. We need to show that 
any pre-image 
$\hat{m}\in\wh{\Phi}^{-1}(L\t^*)$ is fixed under $t$. 
By definition of the action on $M$, any pre-image satisfies 
$t\cdot \hat{m}=g\cdot \hat{m}$ for some $g\in\Om G$. 
By equivariance and since $T$ acts trivially on $L\t^*$, this 
means $\wh{\Phi}(\hat{m})=g\cdot \wh{\Phi}(\hat{m})$. Since
$\Om G$ acts freely, we conclude $g=e$ and therefore 
$t\cdot\hat{m}=\hat{m}$. 
Viewing $L\t^*\to T$ as a principal $\Om T$-bundle, $\wh{F}\to F$ is 
the pull-back bundle under the map $\Phi_F:\,F\to T$. From the 
constructions, it follows that $\wh{F}$ is the Hamiltonian 
$LT$-manifold associated to $F$, proving (b). 

Now view $\t\to T$ as a $\Lambda$-principal sub-bundle of $L\t^*\to T$. 
We have $L\t^*/\Om_0T=\t$ since every
$\Om_0 T$-orbit in $L\t^*$ passes through a unique point in $\t$. Then
(c) follows since $\wh{F}\to F$ and $\widetilde{F}\to F$ are
pull-back bundles with respect to $\Phi_F:\,F\to T$, and from the fact
that the form $\varpi\in \Om^2(L\g^*) $ vanishes if pulled back to
$\t\subset L\t^*$.
\end{proof}

\begin{remark}
By a similar argument, the fixed point 
set $M^g$ of any group element $g\in G$ is a group valued Hamiltonian
$G_g$-manifold, with the pull-backs of $\om$ and $\Phi$ as 2-form 
and moment map.
\end{remark}
Under the assumptions of Theorem \ref{onal}, we are now going to  
explain the ingredients of the fixed point contributions \eqref{eq:LF}. 
First, we will define a level $2(k+c)$  ``$\Spin_c$'' line bundle  
$\L:=L^2\otimes K^{-1}\to\wh{M}$, where $K^{-1}$ is the 
``anti-canonical line bundle'' $K^{-1}$ for Hamiltonian loop 
group manifolds \cite{me:can}. As we will explain, the restriction 
$\L|_{\wt{F}}$ descends to a $T_{2(k+c)}$-equivariant line bundle 
$\L_F\to F$. The element $t$ acts on $\L_F$ as scalar multiplication 
by some element $\zeta_F(t)$, and we will show how to choose a 
square root.

To carry out the details, we need the symplectic cross-section theorem for 
Hamiltonian $LG$-manifolds, cf. \cite{me:lo}. It is an analog 
of the Guillemin-Sternberg cross-section theorem \cite{gu:sy}
for compact groups. For any vertex $\sig$ of $\Alc$,
let $\Alc_\sig\subset\Alc$ denote the complement of the closed face
opposite $\sig$.  Given an arbitrary open face $\sig$, define
$\Alc_\sig$ to be the intersection of all $\Alc_\tau$ with $\tau$ a
vertex of $\ol{\sig}$. Then the flow-outs
$$ 
{U}_\sig:=LG_\sig\cdot\Alc_\sig\subset {\Oma}^1(S^1,\g)\cong L\g^*
$$  
are smooth, finite dimensional submanifolds, and are slices for points
in $\sig$.  The maps $LG\times_{LG_\sig}{U}_\sig\to L\g^*$ are 
embeddings as open subsets and their images form an open cover.
The cross-section theorem states that the pre-images
${Y}_\sig= \wh{\Phi}^{-1}({U}_\sig)$ are $LG_\sig$-invariant,
finite dimensional symplectic submanifolds, and are Hamiltonian
$\wh{LG}_\sig$-manifolds with the restriction $\wh{\Phi}_\sig$ of
$\wh{\Phi}$ as a moment map. (The central circle 
$\U(1)\subset\wh{LG}_\sig$ acts trivially, with moment map $1$.)  
%
%
If
$\nh{L}\to \wh{M}$ is a pre-quantum line bundle, it restricts to an
$\wh{LG}_\sig$-equivariant pre-quantum line bundle for ${Y}_\sig$.
%
In analogy to the finite dimensional setting, we define 
a ``$\Spin_c$''-line bundle $\ca{L}\to\wh{M}$ as a 
tensor-product of $\nh{L}^2$ with an anti-canonical line bundle 
$\nh{K}^{-1}\to \wh{M}$. The notion of anti-canonical line bundle 
for Hamiltonian loop group manifolds was 
introduced in  \cite{me:can}. $\nh{K}^{-1}$ is a level 
$2c$ line bundle with the property that for each cross-section ${Y}_\sig$, 
there is an $\wh{LG}^{(2c)}_\sig$-equivariant 
isomorphism 
\begin{equation}\label{eq:canonical}
 \nh{K}^{-1}|_{{Y}_\sig}\cong K_\sig^{-1}\otimes 
\C_{2(\rho-\rho_\sig,c)}
\end{equation}
where $K_\sig$ is the canonical line bundle for ${Y}_\sig$. 
Here we are using that since $\gamma_\sig=B_c^\sharp(\rho-\rho_\sig)$ 
is contained in $\sig$, the weight $2(\rho-\rho_\sig,c)$ 
defines a 1-dimensional representation of $\wh{LG}^{(2c)}_\sig$.
The conditions \eqref{eq:canonical} are consistent because 
for $\sig\subset\ol{\tau}$, there is an $LG_\tau$-equivariant 
isomorphism
$$ K_\sig^{-1}|_{Y_\tau}=K_\tau^{-1}\otimes \C_{2(\rho_\tau-\rho_\sig)}.$$

Assume $t\in T^{reg}$ and let $F\subseteq M^t$ be a connected  component 
of the fixed point set. By part (c) of Proposition \ref{prop:resext}, 
the action of $\wh{LG}^{(2(k+c))}$ restricts to an action of 
$T_{2(k+c)}\times \wh{\Lambda}^{(2(k+c))}$, and furthermore by part (a)
the central extension $\wh{\Lambda}^{(2(k+c))}$ is trivial. 
Choosing any trivialization (by choosing lifts generators of $\Lambda$), 
we obtain a $T_{2(k+c)}$-equivariant line 
bundle $\L_F\to F$, by setting $\L_F:=\L|_{\wt{F}}/\Lambda$. 
%
%

Suppose now that $t\in T_{k+c}^{reg}$ and  
let $\zeta_F(t)$ be the eigenvalue for the action on $\L_F$. We show 
how to specify a square root $\zeta_F(t)^{1/2}$.
Let 
$$\zeta_{\wt{F}}(t), \mu_{\wt{F}}(t), \kappa_{\wt{F}}(t)
:\,\wt{F}\to\U(1)$$ 
be the (locally constant) eigenvalues for 
the action of $t$ on $\L|_{\wt{F}},\, L|_{\wt{F}},\, \nh{K}|_{\wt{F}}$,
respectively. Then
$\zeta_{\wt{F}}(t)=\mu_{\wt{F}}(t)^2\kappa_{\wt{F}}(t)^{-1}$.
In order to define the square root of $\zeta_{\wt{F}}(t)$ we need to
define the square root $\kappa_{\wt{F}}(t)^{-1/2}$.
Given a face $\sig$ of $\Alc$ and $w\in\Waff$ let ${Y}_{w\sig}:=g\cdot
{Y}_\sig$, where $g\in N_G(T)\ltimes \Lambda\subset LG$ represents
$w$.  It is a finite dimensional symplectic submanifold, invariant
under the action of $LG_{w\sig}:=\Ad_g(LG_\sig)$. Then 
$$ \wh{\Phi}^{-1}(\t)\subset 
\bigcup_{\sig\subset\Alc}\bigcup_{w\in\Waff}
{Y}_{w\sig}
$$ 
so that the intersections $Y_{w\sig}\cap \wt{F}$ cover $\wt{F}$.
By \eqref{eq:canonical}, if $\kappa_{\wt{F}}^{w\sig}(t)^{-1}$ 
is the eigenvalue for the action on the anti-canonical line bundle for 
$Y_{w\sig}$, 
$$ 
\kappa_{\wt{F}}(t)^{-1}\Big|_{Y_{w\sig}\cap \wt{F}}
=\kappa_{\wt{F}}^{w\sig}(t)^{-1}\ 
t^{2 w(\rho-\rho_\sig)}
$$
where $w(\rho-\rho_\sig)$ is defined using the level
$c$ action of $\Waff$.
As in Section \ref{sec:alter} we can define the square root
of $\kappa_{\wt{F}}^{w\sig}(t)^{-1}$. 

\begin{lemma}
\label{lem:root}
There exists a unique 
locally constant $\U(1)$-valued function 
$\kappa_{\wt{F}}(t)^{-1/2}$
on $\wt{F}$ such that 
\begin{equation}\label{eq:squareroot}
\kappa_{\wt{F}}(t)^{-1/2}
\Big|_{Y_{w\sig}\cap \wt{F}}
=(-1)^{\eps(w,\sig)+l(w)}
\kappa_{\wt{F}}^{w\sig}(t)^{-1/2} e^{2\pi i  \l w(\rho-\rho_\sig),v \r}.
\end{equation}
Here $v\in\t$ is the unique vector in $W\cdot\Alc$ 
with $\exp(v)=t$, $l(w)$ is the length of $w$, 
and $\eps(w,\sig)$ is the number of positive roots 
$\alpha\in\mf{R}_{+,\sig}$ of $G_\sig$ (cf. Section \ref{sigplus})
such that $\l w_1\alpha,v\r<0$, where $w_1\in W$ 
is the image of $w$ under the quotient map $\Waff\to W$.  
Under the action of 
$\xi\in \Lambda$,  
$$ \xi^*\kappa_{\ti F}(t)^{-1/2}=t^{B_c^\flat(\xi)}
\kappa_{\ti F}(t)^{-1/2}.
$$
\end{lemma}

\begin{proof}
Note first of all that the right hand side of \eqref{eq:squareroot} 
is well-defined. Indeed, if $w$ is replaced by $w'$ 
with $w\sig=w'\sig$, 
the factor $e^{2\pi i  \l w(\rho-\rho_\sig),v \r}$ 
does not change because $\rho-\rho_\sig\in B_c^\flat(\sig)$ 
is fixed under the level $c$ action of any element of 
$\Waff$ fixing $\sig$, and $l(w)+\eps(w,\sig)$ changes by an even number.
Given faces $\sig\subset\ol{\tau}$ of the alcove and any
$w\in \Waff$, the symplectic normal bundle of 
${Y}_{w\tau}$ inside ${Y}_{w\sig}$ is 
$T\subset LG$-equivariantly isomorphic to $\g_\sig/\g_\tau$, 
with $T$ acting via the isomorphism $w_1^{-1}:\,T\to T$ 
induced by $w_1$. Using the sign convention from 
Section \ref{sec:alter}, the square root of the eigenvalue for the 
action of $t=\exp(v)$ on $\g_\sig/\t$ is given by 
$(-1)^{\eps(w,\sig)}e^{2\pi i \l\rho_\sig,w_1^{-1}v\r}$, 
and similarly for the action on $\g_\tau/\t$. Therefore, 
$$
(-1)^{\eps(w,\sig)}\kappa_{\ti F}^{w\sig}(\ti{m},t)^{-1/2}
=(-1)^{\eps(w,\tau)} 
\kappa_{\ti F}^{w\tau}(\ti{m},t)^{-1/2}
e^{2\pi i \l w_1(\rho_\tau)-w_1(\rho_\sig),v\r}.
$$
Since $w_1(\rho_\tau)-w_1(\rho_\sig)=
-w(\rho-\rho_\tau)+w(\rho-\rho_\sig)$, 
we have shown that the right hand sides of equation \eqref{eq:squareroot}
patch together to a well-defined locally constant function on $\wt{F}$.

The action of $\xi\in\Lambda\subset \Waff$ amounts to replacing $w$ by $\xi\cdot
w$. This does not change $\eps(w,\tau)$, and changes $l(w)$ by 
an even number. The factor $e^{2\pi i \l w(\rho-\rho_\sig),v\r}$
changes by $t^{B_c^\flat(\xi)}$.
\end{proof}

Using Lemma \ref{lem:root}, we define 
$$ \zeta_{\wt{F}}(t)^{1/2}:=\mu_{\wt{F}}(t)
\,\kappa_{\wt{F}}(t)^{-1/2}.$$ 
Under the action of $\xi\in\Lambda$ it
transforms according to
$\xi^*\zeta_{\wt{F}}(t)^{1/2}=
t^{B_{c+k}^\flat(\xi)}\,\zeta_{\wt{F}}(t)^{1/2}.$
But $t^{B_{c+k}^\flat(\xi)}=1$ since $t\in T_{k+c}$. 
Hence $\zeta_{\wt{F}}(t)^{1/2}$ is actually a constant, 
which defines $\zeta_{{F}}(t)^{1/2}$. 

\begin{remark}
The following special case of the definition will be used in our
applications to Verlinde formulas.  Suppose that $F$ contains a point
$m \in \Phinv(e)$, and let $\ti{m}\in\wt{F}$ be the unique 
point in the zero level set mapping to $m$. Then the tangent space 
$T_m M$ is {\em symplectic}, and the quotient map $\wh{M}\to M$ 
induces a $t$-equivariant isomorphism of symplectic vector spaces, 
$T_{\ti{m}}Y_{\{0\}}\cong T_mM$. Hence, choosing   
any $t$-invariant compatible complex structure on $T_mM$ 
and letting $A(t)\in \on{Aut}_\C(T_mM)$ denote the action of $t$, 
\begin{equation} \label{eq:simplekappa}
 \kappa_{\wt{F}}(t,\ti{m})^{-1/2} = 
{\det}_\C(A(t)^{1/2}). 
\end{equation}
If we know in addition that $t$ acts trivially on the 
fiber $L_{\ti{m}}$ (e.g. if $t$ is in the identity component of 
$LG_{\ti{m}}$), we obtain
$$
\zeta_{F}(t)^{1/2} = 
{\det}_\C(A(t)^{1/2})$$
with no explicit reference to the loop group space. 
\end{remark}

\subsection{Alternative version of the fixed point expressions}
\label{sec:altern}
The expression for the fixed point contribution of 
$t\in T^{\on{reg}}_{k+c}$ simplifies if for some $\sig\subset \Alc$,
$$\Phi(F)\subset W\cdot\exp(\Alc_\sig).$$ 
Let $W_\sig$ be the Weyl group of $G_\sig$, that is, the subgroup of
$W$ fixing $\exp(\sig)\subset T$. The connected components of  
$W\cdot\exp(\Alc_\sig)$ are $W$-translates of 
$W_\sig \cdot\exp(\Alc_\sig)$. Let $w\in W$ be such that 
$w(W_\sig\exp(\Alc_\sig))$ contains $\Phi(F)$.
The $\wh{LG}$-equivariant pre-quantum bundle on $\wh{M}$
restricts to a pre-quantum bundle for the Hamiltonian $T$-action on
$Y_{w\sig}$, where $T$ is embedded in $\wh{LG}$ using \eqref{eq:triv}.

\begin{proposition}\label{prop:alternative} 
The fixed point contribution $\chi(\nu_F,t)$ is related to the fixed
point contribution $\chi(\nu_F^{w\sig},t)$ for 
the Hamiltonian $T$-space $Y^{w\sig}$ (defined using \eqref{fp2} or
\eqref{lco}) by
$$
 \chi(\nu_F,t)=\chi(\nu_F^{w\sig},t) 
\f{\D_\C(\g_\sig/\t,w^{-1}t)}{\D_\C(\g/\t,w^{-1}t)}.
$$
In particular, if $\sig=\{0\}$, we have 
$\chi(\nu_F,t)=\chi(\nu_F^\sig,t)$. 
\end{proposition}

\begin{proof}
The projection map $\wh{M}\to M$ restricts to an equivariant
diffeomorphism $\wt{F}\cap Y_{w\sig}\to F$.  Let 
$\ca{L}_{w\sig}$ be $\Spin_c$-line bundle corresponding to 
$Y_{w\sig}$. Then
$$\ca{L}_{w\sig}|_{\wt{F}}\cong \L_F\otimes \C_{2w(\rho-\rho_\sig)}.$$
The normal bundle of $F$ in
$M$ splits $T$-equivariantly into the normal bundle $\nu_F^{w\sig}$ in
$Y_{w\sig}$ and the constant bundle $\g/\g_\sig$. Using
$\D_\R(\g/\g_\sig,t)\D_\R(\g_\sig/\t,t)=\D_\R(\g/\t,t)$ we obtain,
$$\D_\R(\nu_F,t)=\D_\R(\nu_F^{w\sig},t)\f{\D_\R(\g/\t,t)}
{\D_\R(\g_\sig/\t,t)}.$$
Let $\zeta_F^{w\sig}(t)^{1/2}\in \U(1)$ be the square root 
for the action on $\L_{w\sig}$. We have 
\beq 
\f{\zeta_F(t)^{1/2}}{\D_\R(\nu_F,t)}
&=&
\f{\zeta_F^{w\sig}(t)^{1/2}}{\D_\R(\nu_F^{w\sig},t)}
(-1)^{\eps(w,\sig)+l(w)}\,e^{2\pi i \l w(\rho-\rho_\sig),v\r}
\f{\D_\R(\g_\sig/\t,t)}{\D_\R(\g/\t,t)}\\
&=&
\f{\zeta_F^{w\sig}(t)^{1/2}}{\D_\R(\nu_F^{w\sig},t)}
\f{\D_\C(\g_\sig/\t,w^{-1}t)}{\D_\C(\g/\t,w^{-1}t)}.
\eeq
\end{proof}

\section{Proof of the fixed point formula}
\label{sec:theproof}

Our proof of the fixed point formula proceeds in two stages. 
First, we show that in case $\wh{M}$ admits a global cross-section, 
the formula follows from the ``quantization commutes with reduction 
theorem'' applied to the cross-section. In a second step we 
reduce to this case using the method of symplectic cutting. 
\subsection{An identity for level $k$ characters}

Let $G$ be equipped with inner product $B=B_k$ where 
$k\in (\Z_{>0})^s$. 
We will need a Lemma expressing the restrictions 
of irreducible level $k$ characters of $G$
to the group $T_{k+c}$, in terms of characters of 
central extensions $\wh{G}_\sig$ of $G_\sig$ obtained 
as the pull-back of $\wh{LG}^{(k)}$ by 
the map $G_\sig\cong LG_\sig$ (cf. \eqref{eq:inverse}).
For $\sig\subset\ol{\tau}$ we have embeddings $\wh{G}_\tau\subset 
\wh{G}_\sig$, in particular every $\wh{G}_\sig$ contains 
$\wh{T}$ as a maximal torus. 
Let $\wh{T}=T\times\U(1)$ be the trivialization 
obtained by restricting the trivialization \eqref{eq:triv}.
In terms of the 
corresponding splitting $\wh{\t}^*=\t^*\times \R$, the action of
$W_\sig$ on $\wh{\t}^*$ reads
\begin{equation}\label{eq:affineaction}
w_1\cdot(\mu,\tau)=(w_1\mu+\tau 
B^\flat(\gamma_\sig-w_1\gamma_\sig),\tau)
\end{equation}
and a positive Weyl chamber for $\wh{G}_\sig$ is given by
$$
\wh{\t}^*_{\sig,+}
=(\t^*_{\sig,+}\times\{0\})+\R\cdot(B^\flat(\gamma_\sig),1). 
$$
For any level $k$ weight $\mu\in\Lambda^*_k$, the weight 
$(\mu,1)\in\Lambda^*\times\Z$ is contained in the 
positive Weyl chamber for $\wh{G}_\sig$, and hence parametrizes 
an irreducible representation. Consider the restriction of 
its character $\chi_{\mu,\sig}\in C^\infty(\wh{G}_\sig)$ to 
$$ T_{k+c}\subset  T \subset \wh{T}\subseteq \wh{G}_\sig.$$
We identify $W/W_\sig$ with the set of all $w\in W$ such that 
$w(\t_+)\subset \t_{+,\sig}$. Every element in $W$ can be uniquely  
written in the form $w w_1$ with $w\in W/W_\sig$ and $w_1\in W_\sig$. 
\begin{lemma}
\label{lem:characterrestriction}
For all $t\in T_{k+c}$, and all $\mu\in\Lambda^*_k$, 
$$ \chi_\mu(t)=\sum_{w\in W/W_\sig}
\chi_{\mu,\sig}(w^{-1}t)
\f{\D_\C(\g_\sig/\t,w^{-1} t)}{\D_\C(\g/\t,w^{-1} t)}.$$
\end{lemma}

\begin{proof}
By the Weyl character formula,
$$
\chi_\mu(t)=
\sum_{w\in W/W_\sig}\sum_{w_1\in W_\sig}
(-1)^{l(w\,w_1)}\f{t^{ww_1(\mu+\rho)-\rho}}{\D_\C(\g/\t,t)}
=
\sum_{w\in W/W_\sig}\f{\sum_{w_1\in W_\sig}
(-1)^{l(w_1)} (w^{-1}t)^{w_1(\mu+\rho)-\rho}}{\D_\C(\g/\t,w^{-1}t)}.
$$
Given $t\in W/W_\sig$ let $t_1=w^{-1}t$. We claim that the sum 
over $W_\sig$ is just 
$$
\D_\C(\g_\sig/\t,t_1)
\chi_{\mu,\sig}(t_1)=
\sum_{w_1\in W_\sig}(-1)^{l(w_1)}
t_1^{w_1(\mu+\rho_\sig,1)-(\rho_\sig,0)}.
$$
Indeed, by \eqref{eq:affineaction}
and since $\rho_\sig=\rho-B_c^\flat(\gamma_\sig)$, we have
$$w_1(\mu+\rho_\sig,1)-(\rho_\sig,0)
=(w_1(\mu+\rho)-\rho+B_{k+c}^\flat(\gamma_\sig-w_1\gamma_\sig),1).
$$
But $t_1^{B_{k+c}^\flat(\gamma_\sig-w_1\gamma_\sig)}=1$ since 
$t_1\in T_{k+c}$. Hence  
$ t_1^{w_1(\mu+\rho_\sig,1)-(\rho_\sig,0)}
=t_1^{w_1(\mu+\rho)-\rho},
$
proving the claim.
\end{proof}

\subsection{Proof in case of a global cross-section}
We now explain the proof of Theorem \ref{onal} in the special 
case where $(\wh{M},\wh{\om},\wh{\Phi})$ admits a global 
cross-section. That is, we make the assumption that for 
some face $\sig$ of the alcove, the moment polytope is 
contained in $\Alc_\sig$. As a consequence
$$ \wh{M}={LG}\times_{{LG}_\sig}{Y}_\sig.$$
Using the identification $\wh{LG}_\sig\cong \wh{G}_\sig$, we  
view ${Y}_\sig$ as a Hamiltonian $\wh{G}_\sig$-space. 
Clearly, $\wh{M}_\mu=(Y_\sig)_\mu$ for all $\mu\in\Lambda^*_k$. 
Using Lemma \ref{lem:characterrestriction} 
and the ``quantization commutes with reduction'' 
principle (Theorem \ref{qrc}),  
\beq 
\sum_{\mu\in\Lambda^*_k}
\chi(\wh{M}_\mu)\chi_\mu(t)&=&
\sum_{\mu\in\Lambda^*_k}\chi((Y_\sig)_\mu)\,
\sum_{w\in W/W_\sig}\chi_{\mu,\sig}(w^{-1}t)
\f{\D_\C(\g_\sig/\t,w^{-1}t)}{\D_\C(\g/\t,w^{-1}t)}
\\
&=&\sum_{w\in W/W_\sig}
\chi(Y_\sig,w^{-1}t)
\f{\D_\C(\g_\sig/\t,w^{-1}t)}{\D_\C(\g/\t,w^{-1}t)}.
\eeq
Theorem \ref{onal} now follows by an application of the 
fixed point formula to $Y_\sig$, and using 
Proposition \ref{prop:alternative} 
for the fixed point contributions.

\subsection{Proof in the general case}
Our proof of Theorem \ref{onal} in the general case is an application
of symplectic cutting, reviewed in Appendix \ref{sec:symplcut}.

Denote by $\Psi:\,M\to \Alc$ the composition of the map 
$\Phi:\,M\to G$ 
with the quotient
map $G\to G/\Ad(G)=\Alc$.  For suitable polytopes $Q\subset \t$, the
cut spaces $M_Q$ will be obtained by collapsing the boundary of
$\Psi^{-1}(Q)$ in a certain way. The polytopes $Q$ are defined as
follows.

\begin{figure}[htb]
\setlength{\unitlength}{0.005in}
\begingroup\makeatletter\ifx\SetFigFont\undefined
\def\x#1#2#3#4#5#6#7\relax{\def\x{#1#2#3#4#5#6}}%
\expandafter\x\fmtname xxxxxx\relax \def\y{splain}%
\ifx\x\y   
\gdef\SetFigFont#1#2#3{%
  \ifnum #1<17\tiny\else \ifnum #1<20\small\else
  \ifnum #1<24\normalsize\else \ifnum #1<29\large\else
  \ifnum #1<34\Large\else \ifnum #1<41\LARGE\else
     \huge\fi\fi\fi\fi\fi\fi
  \csname #3\endcsname}%
\else
\gdef\SetFigFont#1#2#3{\begingroup
  \count@#1\relax \ifnum 25<\count@\count@25\fi
  \def\x{\endgroup\@setsize\SetFigFont{#2pt}}%
  \expandafter\x
    \csname \romannumeral\the\count@ pt\expandafter\endcsname
    \csname @\romannumeral\the\count@ pt\endcsname
  \csname #3\endcsname}%
\fi
\fi\endgroup
\begin{picture}(570,375)(0,-10)
\thicklines
\path(530,70)(70,70)
\path(530,335)(70,70)
\path(530,70)(530,335)(530,335)
\thinlines
\path(508,299)(488,332)(121,121)(142,90)
\path(488,331)(510,360)(550,360)
	(570,331)(549,299)
\path(140,50)(120,20)(90,0)
	(50,0)(20,20)(0,50)
	(0,90)(20,120)(50,140)
	(85,140)(120,120)
\path(141,50)(141,89)(510,89)
	(510,50)(141,50)
\path(142,89)(509,298)
\path(510,49)(510,89)(551,89)
	(551,49)(510,49)
\path(550,89)(550,297)(509,297)
	(509,89)(550,89)
\end{picture}
\caption{  The polytopes $Q_\sigma$ for $G_2$. The bold-faced line indicates 
the boundary of the Weyl alcove $\Alc$.}
\label{Qfig}
\end{figure}
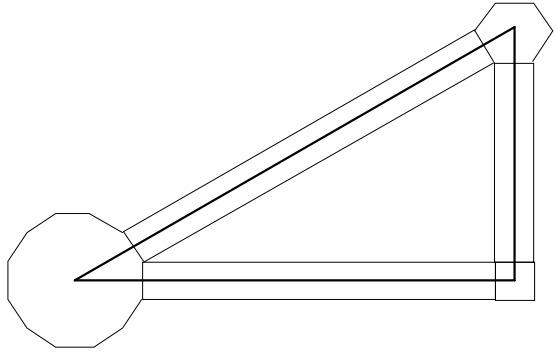

Pick a rational point $\mu\in\Lambda\otimes_\Z\Q$ in the interior of
the alcove, and let $\eps\in\Q$ with $0<\eps<1$.  For any face $\sig$
of $\Alc$, let $Q=Q_\sig$ be the convex hull of all conjugates of
$(1-\eps)\ol{\sig}+\eps\mu\subset\Alc$ under the affine action of
$W_\sig$. (See Figure \ref{Qfig}).  The polytope $Q$ is simplicial;
choose integral labels as in \ref{sec:symplcut1}.  Notice that
$Q\cap\Alc\subset\Alc_\tau$ for all $\tau \subseteq \ol{\sig}$.

The polytope $Q=Q_\sig$ will be called $\Phi$-admissible if it is
$\wh{\Phi}_\sig$-admissible (cf. Appendix \ref{sec:b2}) for $Y_\sig$. 
It is then also $\wh{\Phi}_\tau$-admissible for $Y_\tau$ for each
$\tau\subseteq\ol{\sig}$. The cut spaces satisfy
$(Y_\sig)_Q=G_\sig\times_{G_\tau}(Y_\tau)_Q$, so that the orbifold   
$M_Q:=G\times_{G_\tau}(Y_\tau)_Q$ is independent of the choice
of $\tau$ with $\tau\subseteq \ol{\sig}$. There is a natural map 
$\Psi^{-1}(Q)\to M_Q$ which is a diffeomorphism over 
$\Psi^{-1}(\on{int}(Q))$. 
More generally, for $Q$ a face of $Q_\sig$ we let $(Y_\tau)_Q$ be the
corresponding symplectic sub-orbifold of $(Y_\tau)_{Q_\sig}$, and
$M_Q:=G\times_{G_\tau}(Y_\tau)_Q$ is a sub-orbifold of
$M_{Q_\sig}$. Guided by Lemma \ref{lem:characterrestriction} we
define, for $t\in T_{k+c}$,
$$ \chi(M_Q,t):=\sum_{w\in W/W_\tau} 
\f{\D_\C(\g_\tau /\t,w^{-1} t)}{\D_\C(\g /\t,w^{-1} t)}
\chi((Y_\tau)_Q,w^{-1}t),
$$
which again is independent of the choice of $\tau$.

Let $\ca{Q}$ be the collection 
of all $Q_\sig$, along with their conjugates 
under the action of the affine Weyl group $\Waff$. 
By a generic choice of $\mu,\eps$, we can assume 
that all $Q=Q_\sig$ are $\Phi$-admissible.  
Let $\wt{\ca{Q}}$ be the set of all polytopes 
$Q\in \ca{Q}$, and all of their closed faces. 
The following observation is our starting point for the proof 
of Theorem \ref{onal}. 
\begin{lemma}\label{lem:alternating}
For all $t\in T_{k+c}^{\on{reg}}$, 
\begin{equation}\label{eq:startingpoint}
\chi(M,t)=\sum_{Q\in\wt{\ca{Q}}}(-1)^{\codim Q} \chi(M_Q,t).
\end{equation}
\end{lemma}
\begin{proof}
For all $\mu\in \Lambda^*_k\cap Q$, with $Q\cap\Alc\subset\Alc_\sig$  
we have $((Y_\sig)_Q)_\mu=\wh{M}_\mu$. 
Hence, ``quantization 
commutes with reduction'' (Theorem \ref{qrc})
together with Lemma \ref{lem:characterrestriction} 
shows that 
$$ \chi(M_Q,t)=\sum_{\mu\in\Lambda^*_k\cap Q}\chi(\wh{M}_\mu)
\chi_\mu(t).$$
Let $1_Q$ be the characteristic function of $Q$. Using the Euler
formula 
$$\sum_{Q\in\ti{\ca{Q}} }(-1)^{\codim Q}1_Q(\mu)=1,$$ 
the alternating sum over 
$\chi(M_Q,t)$ equals
$\sum_{\mu\in\Lambda^*_k}\chi(\wh{M}_\mu)\chi_\mu(t)$.
\end{proof}
The orbifold version of the fixed point formula, 
Theorem \ref{th:vergne} in Appendix \ref{sec:orbifold}, expresses 
all indices
$\chi((Y_\sig)_Q,t)$, and therefore all $\chi(M_Q,t)$, as a sum over
fixed point contributions.  Our aim is to identify this sum with the
sum over fixed point contributions $\sum_{F\subseteq M^t}\chi(\nu_F,t)$.
To obtain the required gluing formula, we would like to 
localize further to the fixed point set of the maximal torus $T\subset G$. 
However, a problem arises because the $T_{2(k+c)}$-action on 
$\L_F\to F$ need not extend to a $T$-action. 
Over each $F\cap Y_\sig$ such a $T$-action can be 
introduced by choice of a moment map, however the local $T$-actions 
obtained in this way do not fit together in general. 

In order to get around this problem, we proceed as in \cite[Appendix
A]{me:lo} and consider a {\em second} collection $\ca{S}=\{S\}$ of integral
labeled polytopes in $\t$. The polytopes in $\ca{S}$ are constructed
just like those in $\ca{Q}$, but with $\eps$ replaced by some $\eps'>
\eps$.  By a generic choice of $\mu,\eps,\eps'$ we may assume that all
$S\in\ca{S}$ and all intersections $S\cap Q$ with $S\in\ca{S}$ and
$Q\in\ca{Q}$ are admissible.  Given $S\in\ca{S}$ and $t\in
T_{k+c}^{\on{reg}}$, we define
$$ \chi_S(M,t)=\sum_{F\subseteq M^t}\chi_S(\nu_F,t),$$
where $\chi_S(\nu_F,t)$ is defined by an integral 
similar to $\chi(\nu_F,t)$ (cf. \eqref{eq:LF}), but integrating only over 
the subset $F\cap \Psi^{-1}(S)$: 
$$ \chi_S(\nu_F,t):= 
\zeta_F(t)^{1/2}\int_{F\cap \Psi^{-1}(S)}
\frac{\widehat{A}(F) e^{\hh c_1(\L_F)} }{\D_\R(\nu_F,t)} .
$$
\begin{lemma}\label{lem:corners}
For all $S\in\ca{S}$, the 
integral $\chi_S(\nu_F,t)$ is independent of the choices
of differential form representatives 
$ \widehat{A}(F)$, $c_1(\L_F)$ and $\D_\R(\nu_F,t)$, 
provided these are chosen in such a 
way that for each boundary face $R\subset S$, the pull-back of 
the form to $F\cap \Psi^{-1}(R)$ is $T_R$-basic. We have 
\begin{equation}
\label{eq:sums} 
\chi(\nu_F,t)=\sum_{S\in\ca{S}}\chi_S(\nu_{F},t).
\end{equation}
\end{lemma}
\begin{proof}
The first part follows by observing that the integral can be re-written 
as an integral over the cut space $F_S\subset M_S$ of $F$, that is over the 
image of $F\cap\Psi^{-1}(S)$ in $M_S=\Psi^{-1}(S)/ 
\sim $: 
$$ \chi_S(\nu_F,t)=
\zeta_F(t)^{1/2}\int_{F\cap \Psi^{-1}(S)}
\frac{\widehat{A}(F) e^{\hh c_1(\L_F)} }{\D_\R(\nu_F,t)} 
= 
\zeta_F(t)^{1/2}\int_{F_S}
\frac{\widehat{A}((TF)_S) e^{\hh c_1((\L_F)_S)} }{\D_\R((\nu_F)_S,t)}. 
$$
Here $(TF)_S=TF|_{F\cap \Psi^{-1}(S)}/\sim$ is 
the ``cut'' of $TF$ as explained in Appendix 
\ref{sec:b2}, and similarly for 
$(\nu_F)_S$ and $(\L_F)_S$. Formula \eqref{eq:sums} 
is expressing the integral over $F$ as a sum 
of integrals over all pieces $F\cap \Psi^{-1}(S)$ in the decomposition. 
\end{proof}
The integrals $\chi_S(\nu_F,t)$ can be re-written in terms 
of cross-sections $Y_\sig$. As above, we identify $Y_\sig\subset \wh{M}$ 
with its image under the map $\wh{M}\to M$, and interpret $Y_\sig$ 
as a Hamiltonian $\wh{G}_\sig\cong\wh{LG}_\sig$-space. Let 
$\Psi_\sig:\,Y_\sig\to\Alc$ be the restriction of $\Psi$. 
It can be identified with the composition of 
the moment map $\wh{\Phi}_\sig$ with the projection map 
$\wh{\g}_\sig^*\to\wh{\t}^*_+\supset \Alc\times \{1\}$. 
The intersection $Y_\sig\cap\Psi_\sig^{-1}(S)$ 
is compact, and we can define 
$$\chi_S(Y_\sig,t)=\sum_{F\subseteq (Y_\sig)^t} \chi_S(\nu_F^\sig,t),$$ 
with   
$$ \chi_S(\nu_F^\sig,t)=
\mu_F(t) \int_{F\cap\Psi_\sig^{-1}(S)}\f{\Td(F)
e^{c_1(L_\sig|_F)}}{\D_\C(\nu_F^{\sig},t)},
$$
where $\nu_F^{\sig}$ is the normal bundle of $F$ in $Y_\sig$.
As in Lemma \ref{lem:corners}, the integral does not depend 
on representatives for $\Td(F)$, 
$c_1(L_\sig|_F)$ and $\D_\C(\nu_F^{\sig},t)$, 
provided for each open face $R\subset S$ the pull-backs to
$\Psi_\sig^{-1}(R)\cap F$ descend to the quotient by $T_R$.
Following the argument in Section \ref{sec:altern},
we have
$$ \chi_S(\nu_F,t)=
\sum_{w\in W/W_\sig} 
\f{\D_\C(\g_\sig/\t,w^{-1} t)}{\D_\C(\g /\t,w^{-1} t)}
\chi_S(\nu_F^\sig,w^{-1}t),
$$
hence 
$$ \chi_S(M,t)=
\sum_{w\in W/W_\sig} 
\f{\D_\C(\g_\sig/\t,w^{-1} t)}{\D_\C(\g /\t,w^{-1} t)}
\chi_S(Y_\sig,w^{-1}t).
$$
Over $Y_\sig$, we have a Hamiltonian $T$-action with $T$-equivariant
pre-quantum line bundle. Hence $\chi_S(Y_\sig,w^{-1}t)$ can be written
as a limit of $\chi_S(Y_\sig,w^{-1}(t\exp\xi))$ as $\xi\in\t$
approaches $0$, and by the Berline-Vergne formula the
integral defining $\chi_S(Y_\sig,w^{-1}(t\exp\xi))$ localizes to the
fixed point set of $T$.  The details of this approach are given in
Appendix \ref{sec:symplcut}.  In particular Proposition
\ref{SemiLocalGluing} allows us to re-write $\chi_S(Y_\sig,w^{-1}t)$
as an alternating sum over the corresponding terms for the cut spaces
$(Y_\sig)_Q$:
\begin{equation}\label{eq:inapp}
\chi_S(Y_\sigma, w^{-1}t)=
\sum_{Q\in\ti{\ca{Q}}}
(-1)^{\codim Q} 
\chi_S((Y_\sigma)_Q, w^{-1}t).
\end{equation}
From \eqref{eq:inapp} we obtain,
\beq
\lefteqn{
\sum_{F\subseteq M^t}\chi(\nu_F,t)=\sum_{S\in\ca{S}}\chi_S(M,t)}
\\
&=&
\sum_{S\in\ca{S}}\sum_{Q\in\ti{\ca{Q}}}(-1)^{\codim Q} 
\Big(\sum_{w\in W/W_\sig}
\f{\D_\C(\g_\sig/\t,w^{-1} t)}{\D_\C(\g/\t,w^{-1}t)}
\chi_S((Y_\sig)_Q,w^{-1}t)\Big)
\\
&=&
\sum_{Q\in\ti{\ca{Q}}}
(-1)^{\codim Q} 
\Big(\sum_{w\in W/W_\sig}
\f{\D_\C(\g_\sig/\t,w^{-1} t)}{\D_\C(\g /\t,w^{-1} t)}
\chi((Y_\sig)_Q,w^{-1}t)\Big)\\
&=&\sum_{Q\in\ti{\ca{Q}}}
(-1)^{\codim Q} \chi(M_Q,t)\\&=&
\chi(M,t).
\eeq
This completes the proof of Theorem \ref{onal}.

\begin{appendix}
\section{The equivariant index theorem for orbifolds}
\label{sec:orbifold}
The equivariant version of Kawasaki's index theorem for orbifolds 
is due to  M. Vergne \cite{ve:eq}. A good reference 
is Chapter 14 in Duistermaat's book \cite{du:he}; more 
information can be found in \cite[Section 3]{me:sym}.
We follow the conventions for the definition of a $G$-orbifold 
$M$ as given in \cite{du:he}. 

The Kawasaki-Vergne formula expresses the equivariant index as an
integral over connected components of a certain orbifold
$\ti{M}^g$. There is a natural surjection from $\ti{M}^g$ onto the
fixed point set $M^g$ of $g$, the latter however is not in general a
sub-orbifold of $M$:
\begin{example}\label{ex:orbifold}
Let $G=S^1$ act on $\C^2$ by $e^{i\phi}\cdot(z_1,z_2)=
(e^{i\phi}z_1,e^{2 i\phi}z_2)$, and let $\Z_2$ act by  
$(z_1,z_2)\mapsto -(z_1,z_2)$. The $G$-action descends to $M=\C^2/\Z_2$. 
The fixed point set of $g=e^{i\pi}\in G$ is 
$M^g=\{(z_1,z_2)|\,z_1z_2=0\}/\Z_2$, which is not a 
sub-orbifold of $M$.  
\end{example}
Given $m\in M^g$, let $(V,\Gamma,p)$ be a local orbifold chart around 
$m$. Thus $V$ is an open subset of $\R^n$, $\Gamma$ a finite group acting 
on $V$, and $p:\, V/\Gamma\to M$ a homeomorphism onto an open neighborhood 
of $m$. The action of $g$ on $V/\Gamma$ corresponds to some action 
on $V$, together with an automorphism $\phi$ of $\Gamma$ such that 
such that $\gamma\cdot g\cdot x = g\cdot \phi(\gamma)\cdot x$
for all $x\in V,\ \gamma\in\Gamma$. 
Let $\Gamma$ act on 
$$ \ti{V}^g:=\coprod_{\gamma\in\Gamma} 
V^{\gamma g}\times \{\gamma\}$$ 
by $\gamma_1\cdot(x,\gamma)=(\gamma_1 x,\ \gamma_1\gamma\phi(\gamma_1)^{-1})$. 
The orbifold $\ti{M}^g$ is obtained by gluing together 
the orbifolds $\wt{V/\Gamma}^g:=\ti{V}^g/\Gamma$. Usually it has 
a number of connected components of different dimensions. The natural 
maps $\ti{V}^g\to V$ descend to surjective maps $\wt{V/\Gamma}^g\to 
(V/\Gamma)^g$, which patch together to a surjective map 
$\ti{M}^g\to M^g$. 
In our applications, the groups $\Gamma$ are 
abelian and the automorphism $\phi$ is trivial. 
%
%
In Example \ref{ex:orbifold}, $\ti{M}^g$ has two connected 
components:
$$\ti{M}^g =\{(z_1,z_2)|z_2=0\}/\Z_2 \sqcup
\{(z_1,z_2)|z_1=0\}/\Z_2.$$
%
%
%
%
%
%
Suppose now that $(M,\om,\Phi)$ is a Hamiltonian 
$G$-orbifold. Then all connected components of 
$\ti{M}^g$ are symplectic manifolds, with symplectic form 
the pull-back of $\om$ under the map $\ti{M}^g\to M$. 
If $M$ carries a $G$-equivariant pre-quantum line bundle $L\to M$ 
then its pull-back $\ti{L}\to \ti{M}^g$ is a pre-quantum line bundle for 
this symplectic structure. Local charts for $\ti{L}$ are obtained 
from orbifold charts $(V,\Gamma,p)$ for $M$. In any such 
chart, $L$ is given by a 
$\Gamma$-equivariant pre-quantum line bundle $L_V\to 
V$, and the pull-back $L_{\ti{V}^g}$ to $\ti{V}^g$ is a 
$\Gamma$-equivariant line bundle, defining local charts for 
$\ti{L}$.  

At any point $(v,\gamma)\in\ti{V}^g$,\,  $\gamma g$ acts 
on the fiber over $(v,\gamma)$. The weight for this action 
depends only on the connected 
component $\wt{F}$ of $\ti{M}^g$ containing $(v,\gamma)$, 
and will be denoted $\mu_{\wt{F}}(g)\in \U(1)$. 

For any connected component $\wt{F}$, let $\d_{\wt{F}}$ denote its
{\em multiplicity}, that is the number of elements in the orbifold
isotropy group for a point in its smooth part. Let $\nu_{\wt{F}}\to\wt{F}$ be
the normal bundle for the immersion $\wt{F}\to M$. 
In local orbifold charts
$\ti{V}^g$, it is given as the normal bundle $\nu_{\ti{V}^g}$ of
$\ti{V}^g$ in $V\times\Gamma$. Again, its fiber at $(v,\gamma)$
carries an action of $\gamma g$. We let $\ti{\D}_\C(\nu_{\wt{F}},g)\in
\Om(\wt{F})$ be the differential form given in local charts by
$\D_\C(\nu_{\ti{V}^g},\gamma\,g)$, using the definition of $\D_\C$
given in Section \ref{sec:fixed}.  Finally, we can state the fixed
point theorem for this particular case:

\begin{theorem}[Vergne]\label{th:vergne}
Let $(M,\om,\Phi)$ be a pre-quantized Hamiltonian $G$-orbifold. 
For any $g\in G$, the following fixed point formula holds: 
\begin{equation}
\chi(M,g)=\sum_{\wt{F}\subseteq \ti{M}^g} \chi(\nu_{\wt{F}},g)
\end{equation}
where 
\begin{equation}
\chi(\nu_{\wt{F}},g)=\f{1}{d_{\wt{F}}}
\int_{\wt{F}}\f{
\Td(\wt{F})\Ch(\ti{L})}{\ti{\D}_\C(\nu_{\wt{F}},g)}
\mu_{\wt{F}}(g).
\end{equation}
\end{theorem}

We will also need a slightly more general version, expressing the 
index $\chi(M,g\,e^\xi)$ where $\xi$ is a sufficiently small 
element in the Lie algebra of the centralizer of $g$. Let 
$\Ch(\ti{L},\xi)$ and $\Td(\wt{F},\xi)$ and $\ti{\D}_\C(\nu_{\wt{F}},g,\xi)$ 
denote the equivariant extensions, defined by replacing curvatures 
by equivariant curvatures in the definitions. (See e.g. the 
book \cite{be:he}). Then Vergne in \cite{ve:eq} proves the 
more general formula $\chi(M,g\,e^\xi)=
\sum_{\wt{F}\subseteq \ti{M}^g} \chi(\nu_{\wt{F}},g,\xi)$ with 
\begin{equation}
\chi(\nu_{\wt{F}},g,\xi)=\f{1}{d_{\wt{F}}}\ 
\int_{\wt{F}}
\f{\Td(\wt{F},\xi)\Ch(\ti{L},\xi)}
{\ti{\D}_\C(\nu_{\wt{F}},g,\xi)}
\mu_{\wt{F}}(g).
\end{equation}

\section{Symplectic Cutting}\label{sec:symplcut} 
In this Section we explain the non-abelian version of 
Lerman's technique of symplectic cutting. In a nutshell, 
the method associates to any compact Hamiltonian $G$-space 
$(M,\om,\Phi)$, and certain ``$\Phi$-admissible'' polytopes
$Q\subset \t^*$, a Hamiltonian $G$-orbifold $(M_Q,\om_Q,\Phi_Q)$
with moment polytope $\Phi_Q(M_Q)\cap\t^*_+=
\Phi(M)\cap\t^*_+\cap Q$. 
The space $M_Q$ is obtained from $G\cdot\Phi^{-1}(Q\cap
\t^*_+)\subset M$ 
by collapsing the boundary in a certain way.

\subsection{Labeled Polytopes}
\label{sec:symplcut1}
Let $T$ be a torus, with lattice $\Lambda\subset\t$.  A (rational)
{\em polyhedron} $Q\subset \t^*$ is a finite intersection of 
half spaces
$$ Q=\bigcap_{j=1}^N\  \{\mu\in\t^*|\,\,\l\mu,v_j\r\ge r_j\},$$
where $N$ is the number of codimension 1 faces, $v_j\in\Lambda$ are
non-zero lattice vectors and $r_j\in\R$. Compact polyhedra will be
called polytopes.  $Q$ is called {\em simplicial} if for all $\mu\in
Q$, the vectors $v_j$ for which $\l\mu,v_j\r=r_j$ are linearly
independent.  Following \cite{le:ha}, we define a {\em labeled}
polyhedron to be a polyhedron $Q\subset\t^*$ with a choice of inward
pointing normal vectors $v_j\in\Lambda$ for each codimension 1
face. From $Q$ and the labels $v_j$ one recovers the defining
inequalities $\l\mu,v_j\r\ge r_j$. We call a labeled polyhedron {\em
integral} if all $r_j\in\Z$.  Note this does not imply that the
vertices of $Q$ are integral.

Associated to any codimension $k$ face $S$ of a simplicial labeled
polyhedron $Q$ is a $k$-dimensional sub-torus $T_S\subset T$,
with Lie algebra $\t_S$ the space orthogonal to $S$.  Letting
$(v_{i_1,},\ldots,v_{i_k})$ be the labels of codimension 1 faces
containing $S$, the map $\R^k\to \t_S$ which takes the $j$th standard
basis vector to $v_{i_j}$ defines a covering $(S^1)^k\to T_S$, the
kernel of which is isomorphic to the quotient of $\Lambda\cap\t_S$ by
the lattice generated by the vectors $v_{i_j}$.

\subsection{Non-abelian cutting}
\label{sec:b2}
Let $(M,\om,\Phi)$ be a connected Hamiltonian $G$-manifold. 
We denote by $\Psi:\,M\to \t^*_+$ the composition of 
the moment map $\Phi$ with the quotient map 
$\g^*\to \g^*/\Ad(G)=\t^*_+$. It is well-known that over
$\Psi^{-1}(\on{int}(\t^*_+))$, the map $\Psi$ is smooth 
and generates a $G$-equivariant Hamiltonian $T$-action. 
More generally, if $\sig$ is an open face of $\t^*_+$ and 
$Z(G_\sig)\subseteq T$ the center of its centralizer, the composition 
of $\Psi$ with projection to $\mf{z}(G_\sig)^*$ is smooth near 
$\Psi^{-1}(\sig)$, and generates a $G$-equivariant Hamiltonian 
$Z(G_\sig)$-action. 

Let $Q\subset\t^*$ be a simplicial labeled  polyhedron, with the
property that $\Psi^{-1}(Q)$ is compact and connected. Suppose 
that 
\begin{equation}\label{Property1}
\mbox{If $S\subset Q$ and $\sig\subset\t^*_+$ are open faces 
with $S\cap\sig\cap \Psi(M)\not=\emptyset$, then 
$T_S\subseteq Z(G_\sig)$.}
\end{equation}
It then follows that on a $G$-invariant 
neighborhood of $\Psi^{-1}(S)$, the composition 
of $\Psi$ with projection $p_S:\,\t^*\to\t^*_S$ generates 
a Hamiltonian $T_S$-action. 
The polyhedron $Q$ will be called $\Phi$-admissible if it 
has the property \eqref{Property1}, and in addition satisfies
\begin{equation}\label{Property2}
\mbox
{The action of $T_S$ on $\Psi^{-1}(S)$ is locally free.}
\end{equation}
Given a polyhedron $Q$ satisfying condition \eqref{Property1}, 
condition \eqref{Property2} 
can be achieved by an arbitrarily small perturbation 
of the parameters $r_j$. Assuming \eqref{Property1}, 
\eqref{Property2},  
choose a $G$-invariant neighborhood $U_S$ of $\Psi^{-1}(S)$ 
on which the action of $T_S$ is locally free.

Let $(S^1)^k$ act on $U_S$ by means of the covering, 
$(S^1)^k\to T_S$. As a moment map $\phi_S$ for this action we take 
the moment map $p_S\circ \Psi$ for the $T_S$-action, shifted 
by $\nu_S=p_S(S)$: 
$$\phi_S:=p_S\circ \Psi-\nu_S:\,U_S\to \t_S^*\cong (\R^k)^*.$$
The various torus actions and moment maps are compatible, 
in the sense that if $S_1\subset \ol{S_2}$, the restriction of 
$\phi_{S_2}$ to $U_{S_1}\cap U_{S_2}$ is a  
component of $\phi_{S_1}$. 
For any face $S$ consider the symplectic quotients
$$ (U_S)_Q:=(U_S\times \C^k) \qu (S^1)^k$$ 
under the diagonal action, using the standard action on
$\C^k$. From the canonical isomorphism,  
$$ \Psi^{-1}(Q)\cap U_S=(U_S\times \ti{\C}^k) \qu (S^1)^k$$ 
where $\ti{\C}$ is the symplectic manifold with boundary 
$\ti{\C}=S^1\times\R_+\subset T^*(S^1)$, one sees that 
there is a canonical surjective map 
$\Psi^{-1}(Q)\cap U_S \to (U_S)_Q$, which is a symplectomorphism
over $\Psi^{-1}(\on{int}(Q))$.
 
One obtains a Hamiltonian 
$G$-orbifold $(M_Q,\om_Q,\Phi_Q)$, called the {\em cut space}, 
by gluing the open subsets 
$(U_S)_Q$. The cut space is a union of $G$-invariant 
symplectic sub-orbifolds $M_S:=\Psi^{-1}(S)/T_S$, and 
there is a natural surjection $\Psi^{-1}(Q)\to M_Q$ 
which restricts to the quotient maps $\Psi^{-1}(S)\to M_S$.
The normal bundle 
of $M_S$ in $M_Q$ is the associated orbifold bundle,  
$$\nu_S^Q=\Psi^{-1}(S)\times_{T_S} (\C^k/\Gamma_S) \to M_S,$$
where $\Gamma_S$ is the kernel of the homomorphism $(S^1)^k\to T_S$, 
and the action of $T_S$ is induced from the natural 
$(S^1)^k$-action on $\C^k$.

We now extend the cutting construction to $G$-equivariant vector
bundles $E \to M$.  Suppose that for all open faces $S\subset Q$ the
$T_S$-action on $U_S$ lifts to a $G$-equivariant action of its cover
$(S^1)^k\to T_S$, and that for faces $S_1,S_2$ with $S_1\subset
\ol{S_2}$, these torus actions are compatible in the natural
way. Define $(E|U_S)_Q\to (U_S)_Q$ by pulling $E$ back to $U_S\times
\C^k$, restricting to the zero level set for the $(S^1)^k$-action, and
taking the quotient. These local bundles glue together to give a
$G$-equivariant orbifold bundle $E_Q\to M_Q$. Its restriction to $M_S$
is $E_S:=(E|\Psi^{-1}(S))/(S^1)^k$.

We note that the cut $(TM)_Q$ of the tangent bundle $TM$ is not
isomorphic to $T(M_Q)$. Indeed, for any face $S$ of $Q$, 
\begin{equation}\label{CutTangentBundles}
(TM)_Q|_{M_S}=T(M_S)\oplus\C^k\ \ \ \ 
\mbox{ while }\ \ \ \ 
T(M_Q)|_{M_S}=T(M_S)\oplus \nu_S^Q.
\end{equation}
Suppose $L\to M$ is a $G$-equivariant pre-quantum line bundle. 
The $T_S$-action on $U_S$ admits a $G$-equivariant pre-quantum 
lift (with respect to the moment map $p_S\circ \Psi$)
to $L|U_S$. Since the moment map for $(S^1)^k$ is obtained by 
shift by $\nu_S$, the $(S^1)^k$-action admits a pre-quantum lift 
if and only if $\nu_S\in\t_S^*\cong\R^k$ 
is contained in the lattice $\Z^k$. Clearly, this is the case if 
$Q$ is an {\em integral} labeled polyhedron. Since the trivial 
line bundle is pre-quantum for $\C^k$, the cut bundle 
$L_Q\to M_Q$ becomes then a $G$-equivariant pre-quantum bundle
for the cut space. 

\subsection{The gluing formula}
Let $(M,\om,\Phi)$ be a compact, connected 
Hamiltonian $G$-manifold with pre-quantum 
line bundle $L\to M$. A $\Phi$-{\em admissible labeled polyhedral 
subdivision} of $\t^*$ is a collection of $\Phi$-admissible labeled 
polyhedra $\ca{Q}=\{Q\}$ such that the collection covers $\t^*$, the 
intersection of any two polyhedra is either empty or is a face 
of each, and the labels attached to a common codimension 1 face
of any two polyhedra coincide up to sign. We call $\ca{Q}$ 
integral if all of the polytopes $Q\in\ca{Q}$ are integral 
labeled polyhedra. 
Let $\wt{\ca{Q}}$ be the collection of all closed faces of 
$Q\in \ca{Q}$; thus 
$\ca{Q}$ are the top-dimensional polyhedra in $\wt{\ca{Q}}$.
To every $Q\in\wt{\ca{Q}}$ corresponds a pre-quantized Hamiltonian
$G$-orbifold $M_Q$ which is a symplectic sub-orbifold in 
each $M_{Q'}$, such that $Q$ is a closed face of 
$Q'\in \wt{\ca{Q}}$.  
One has the following gluing formula \cite{me:sym} 
relating the $\Spin_c$-indices 
of the cut spaces:
\begin{equation}\label{Gluing}
\chi(M,g)=\sum_{Q\in\wt{\ca{Q}}}(-1)^{\codim Q}\chi(M_Q,g).
\end{equation}
In this paper we  need a refined version
of \eqref{Gluing}, along ideas developed in
\cite{me:lo}. Suppose $(M,\om,\Phi)$ is a pre-quantized 
Hamiltonian $G$-manifold, and $S$ a $\Phi$-admissible
integral labeled polyhedron. Consider the expression
\begin{equation}\label{IntegralOverS}
\chi_S(M,g)=\sum_{F\subseteq M^g} 
\mu_F(g) \int_{F\cap \Psi^{-1}(S)}\f{\Td(F)\Ch(L)}
{\D_\C(\nu_F,g)},
\end{equation}
where the representatives for the Todd class, Chern class 
and $\D_\C(\nu_F,g)$ are chosen in such a way that for all 
faces $R\subset S$, the pull-back to the submanifold $\Psi^{-1}(R)$ 
descends to a form on $M_R=\Psi^{-1}(R)/T_R$. Provided that the 
representatives satisfy this boundary condition,  
\eqref{IntegralOverS} is independent of their choice because the integral 
can be re-written as an integral over the cut space 
$F_S\subset (M^t)_S\subseteq M_S^t$: 
\begin{equation}\label{eq:intovercut} 
 \chi_S(M,g)=\sum_{F\subseteq M^g}\mu_F(g)\int_{F_S}\f{\Td((TF)_S)
\Ch(L_S)}{\D_\C((\nu_F)_S,g)}.
\end{equation}
Suppose that $\ca{Q}=\{Q\}$ is a $\Phi$-admissible, integral, polyhedral
subdivision, and also that all intersections $S\cap Q$ 
with $Q\in\ca{Q}$ are admissible. 
For any $Q\in\ti{Q}$ let $\chi_S(M_Q,g)$ be defined by a formula  
similar to \eqref{IntegralOverS}, by taking the 
integral in the Kawasaki-Vergne formula only over the part $\wt{F}_S$ 
mapping to $S$. 
\begin{equation}\label{eq:intovercut1} 
 \chi_S(M_Q,g)=\sum_{F\subseteq (M_Q)^g}
\mu_{\wt{F}}(g)\int_{\wt{F}_S}\f{\Td((T\wt{F})_S)
\Ch((\ti{L}_Q)_S)}{\ti{\D}_\C((\nu_{\wt{F}})_S,g)}.
\end{equation}
The following Proposition extends 
formula (34) in \cite{me:lo}
to the equivariant case:
\begin{proposition}\label{SemiLocalGluing}
One has the gluing formula, 
\begin{equation} 
\chi_S(M,g)=\sum_{Q\in\wt{\ca{Q}}}(-1)^{\codim Q}\chi_S(M_Q,g).
\end{equation}
\end{proposition}

\begin{proof}
The proof is an extension of the argument 
given in \cite[p.465]{me:lo}. We may assume $g=t\in T$. 
Observe that all of the characteristic forms in 
\eqref{eq:intovercut} 
admit $T$-equivariant extensions. Hence we can write 
$\chi_S(M,t)=\lim_{\xi\to 0}\chi_S(M,t,\xi)$ where 
$$ \chi_S(M,t,\xi)=\sum_{F\subseteq M^t}\mu_F(t)
\int_{F_S}\f{\Td((TF)_S,\xi)
\Ch(L_S,\xi)}{\D_\C((\nu_F)_S,t,\xi)}.
$$
Let us apply the Berline-Vergne localization formula for orbifolds 
(cf. \cite{me:sym}) to this expression. Let $X$ be a fixed point 
manifold for the $T$-action on $(M^t)_S$. 
Letting $\Psi_S:\,M_S\to \t^*_+$ be the map induced by $\Psi$, 
it follows that $\Psi_S(X)$ is a point. 
Let $R\subset S$ be the unique open face containing $\Psi_S(X)$, 
and let $F\subseteq M^t$ be the unique connected component with 
$X\subseteq F_S$. Recall that by \eqref{CutTangentBundles}
the restriction of $(TM)_S$ to $M_R$ is the tangent 
bundle of $M_R$ plus a trivial bundle. Similarly the restriction of 
$(TF)_S$ is the tangent bundle to $X$ plus a trivial bundle, 
and the restriction of $(\nu_F)_S$ is the normal bundle $\nu_X^R$ 
to $X$ in $M_R$. On the other hand, the normal bundle 
of $X$ in $F_S$ is the pull-back of the normal bundle 
$\nu_R^S$ of $M_R$ in $M_S$. We therefore obtain the formula
$$ 
\chi_S(M,t,\xi)=
\sum_{X} \mu_{X}(t\exp\xi) 
\f{1}{d_X} \int_X 
\f{ \Td(X) \Ch(L_R)}
{\D_\C(\nu_X^R,t\exp\xi) \Eul(\nu_R^S,\xi)}.
$$
We obtain similar formulas for all of the cut spaces $M_Q$:
$$ 
\chi_S(M_Q,t,\xi)=
\sum_{X} \mu_{X}(t\exp\xi) 
\f{1}{d_{X}} \int_{X} 
\f{\Td(X) \Ch(\ti{L}_R)    }
{\ti{\D}_\C(\nu_{X}^R,t\exp\xi) \Eul(\nu_R^S,\xi)        },
$$
where the sum is over all connected components $X$ of the $T$-fixed point set 
of $\wt{F}_S$, for all $F\subseteq M_Q^t$. 
If $X$ is a fixed point component for the $T$-action on $(M^t)_S$, 
then $\Psi(X)$ is contained in the interior of a unique 
top-dimensional polyhedron $Q\in\ca{Q}$. Hence, the fixed point 
contribution appears exactly once as a $T$-fixed point 
contribution of the sum 
$\sum_{Q\in\ca{Q}}(-1)^{\codim Q}\chi_S(M_Q,t,\xi)$.
We must show that the remaining fixed point contributions cancel.
These other integrals are over connected components $X$ of $T$-fixed
point sets of $\wt{F}_S$, where $\wt{F}\subseteq \ti{M}_Q^t$ is
defined as in Appendix \ref{sec:orbifold}.  These fixed point
components can be organized as follows. Consider the finite subset of
points of the form $(\Psi_Q)_S(X)\in\Alc$. Given any such point, there
exists a unique open face $R$ of $S$ containing it. If
$R=\on{int}(S)$, then $X$ is simply one of the $T$-fixed point
orbifolds for $M_Q$, and the cancellation of the corresponding fixed
point contributions is just the gluing formula \cite[Theorem
5.4]{me:sym}. In case $R$ is a proper face of $S$, the argument from
\cite{me:sym} carries over without essential change, the reason being
that the fixed point contributions look exactly like fixed point
contributions for $M_Q$, except for the extra factor 
$\Eul(\nu_R^S,\xi)^{-1}$ which appears in all of these integrals.
\end{proof}

\end{appendix}

\end{document}